\documentclass[a4paper,12pt]{amsart}
\usepackage{amssymb}
\usepackage{amsmath}
\usepackage{hyperref}
\usepackage{bbm}
\usepackage{mathrsfs}
\usepackage{mathtools}
\usepackage{bm}
\usepackage{geometry}
\usepackage{colordvi}
\usepackage{color}
\usepackage{tikz}

\usepackage{tikz-cd}
\allowdisplaybreaks

\numberwithin{equation}{section}
\newcommand{\fghopf}[1]{(FinHopf)}
\theoremstyle{definition}
\newtheorem{theorem}{Theorem}[section]
\newtheorem{corollary}[theorem]{Corollary}
\newtheorem{lemma}[theorem]{Lemma}
\newtheorem{definition}[theorem]{Definition}
\newtheorem{proposition}[theorem]{Proposition}
\newtheorem{example}[theorem]{Example}
\newtheorem{remark}[theorem]{Remark}
\newtheorem{que}[theorem]{Question}

\begin{document}


\title{Chevalley property and discriminant ideals of Cayley-Hamilton Hopf Algebras}



\author{Yimin Huang}
\address{School of Mathematical Sciences, Fudan University, Shanghai 200433, China}
\email{21110180008@m.fudan.edu.cn}

\author{Zhongkai Mi}
\address{Shanghai Center for Mathematical Sciences, Fudan University, Shanghai 200438, China}
\email{zhongkai{\_}mi@fudan.edu.cn}

\author{Tiancheng Qi}
\address{School of Mathematical Sciences, Fudan University, Shanghai 200433, China}
\email{tcqi21@m.fudan.edu.cn}

\author{Quanshui Wu}
\address{School of Mathematical Sciences, Fudan University, Shanghai 200433, China}
\email{qswu@fudan.edu.cn}

\begin{abstract}
For any affine Hopf algebra $H$ which admits a large central Hopf subalgebra,  $H$ can be endowed with a Cayley-Hamilton Hopf algebra structure in the sense of De Concini-Procesi-Reshetikhin-Rosso. The category of finite-dimensional modules over any fiber algebra of $H$ is proved to be an indecomposable exact module category over the tensor category of finite-dimensional modules over the identity fiber algebra $H/\mathfrak{m}_{\overline{\varepsilon}}H$ of $H$. 
For any affine Cayley-Hamilton Hopf algebra $(H,C,\text{tr})$ such that  $H/\mathfrak{m}_{\overline{\varepsilon}}H$ has the Chevalley property, it is proved that  if the zero locus of a discriminant ideal of $(H,C,\text{tr})$ is non-empty then it contains the orbit of the identity element of the affine algebraic group  $\text{maxSpec}C$ under the left (or right) winding automorphism group action. Its proof relies on the fact that $H/\mathfrak{m}_{\overline{\varepsilon}}H$ has the Chevalley property if and only if  the $\overline{\varepsilon}$-Chevalley locus of $(H,C)$ coincides with $\text{maxSpec}C$.



Then, we provide a description of the zero locus of the lowest discriminant ideal of $(H,C,\text{tr})$. It is proved that the lowest discriminant ideal of $(H,C,\text{tr})$ is of level $\text{FPdim}(\text{Gr}(H/\mathfrak{m}_{\overline{\varepsilon}}H))+1$, where $\text{Gr}(H/\mathfrak{m}_{\overline{\varepsilon}}H)$ is the Grothendieck ring of the finite-dimensional Hopf algebra $H/\mathfrak{m}_{\overline{\varepsilon}}H$  and $\text{FPdim}(\text{Gr}(H/\mathfrak{m}_{\overline{\varepsilon}}H))$ is the Frobenius-Perron dimension of $\text{Gr}(H/\mathfrak{m}_{\overline{\varepsilon}}H)$. Some recent results of Mi-Wu-Yakimov about lowest discriminant ideals are generalized. We also prove that all the discriminant ideals are trivial if $H$ has the Chevalley property.

\end{abstract}


\thanks{This research has been supported by the NSFC (Grant No. 12471032), the National Key Research and Development Program of China (Grant No. 2020YFA0713200) and the China Postdoctoral Science Foundation (Grant No. 2025M773084 to Z.M.).}

\maketitle


\section*{Introduction} 
A \textit{module-finite algebra} is an algebra $A$ 
with a central subalgebra $C$ such that $A$ is a finitely generated module over $C$. 
For any \textit{affine} (finitely generated as an algebra) module-finite algebra $A$ over $C$, each $\mathfrak{m}\in \text{maxSpec}C$ (the maximal spectrum of $C$) corresponds to a finite-dimensional algebra $A/\mathfrak{m}A$, which is called the \textit{fiber algebra} of $A$ at the maximal ideal $\mathfrak{m}$. 
It is a standard fact that the set of irreducible representations of $A$ is the disjoint union of  the irreducible representations over its fiber algebras. 

Module-finite algebras arise naturally in various fields of mathematics, say, the De Concini-Kac quantized enveloping algebras $\mathcal{U}_{\epsilon}(\mathfrak{g})$ of finite-dimensional complex semisimple Lie algebras $\mathfrak{g}$ at roots of unity $\epsilon$ \cite{MR1103601,MR1288995} and the quantized coordinate rings $\mathcal{O}_{\epsilon}(G)$ of simply connected complex semisimple algebraic groups $G$ at roots of unity $\epsilon$ \cite{MR1296515} in quantum group theory;  non-commutative crepant resolutions of affine Gorenstein normal domains \cite {MR2077594} in non-commutative resolution theory; Kauffman bracket skein algebras $\mathscr{S}_{\epsilon}(\Sigma,\mathcal{P})$ of marked surfaces $(\Sigma,\mathcal{P})$ at roots of unity $\epsilon$ \cite{MR4264235} in knot theory, and symplectic reflection algebras $H_{t,\boldsymbol{c}}$ with parameter $t=0$ \cite{MR1881922} in symplectic resolution theory.

A $C$-linear map $\text{tr}: A\to C$ is called a \textit{trace map} from $A$ to $C$ if $\text{tr}(ab)=\text{tr}(ba)$ for all $a,b\in A$. The triple $(A,C,\text{tr})$ is called an \textit{algebra with trace} \cite{MR883869}. In fact, each of the above classes of algebras, when equipped with suitable trace maps, falls within the framework of \textit{Cayley-Hamilton algebras} in the sense of Procesi \cite{MR883869,MR1288995}, see \S \ref{sec-CH algs} for details. Cayley-Hamilton algebras provide a unified approach to studying their representation theory and ring-theoretic properties, see for example \cite{MR1288995,MR2178656,MR3886192,MR4334165}.

Over the past three decades, module-finite Hopf algebras have attracted much interest \cite{MR4600057,MR4201485,MR1482982,MR1863398,MR1676211,MR2661247,MR1296515,MR1288995,MR2178656,
mi2025lowest}.
Among them, many important module-finite Hopf algebras are affine and admit \textit{large} central Hopf subalgebras (that is, they are finitely generated modules over some central Hopf subalgebras). Most quantum groups at roots of unity fall into this class of Hopf algebras \cite{MR1103601,MR1288995,MR1296515,MR2178656,MR4600057}.  




The first result in this paper reveals that all affine Hopf algebras that admit large central Hopf subalgebras possess a Cayley-Hamilton algebra structure.



\begin{theorem}
[Proposition \ref{prop-equfibalg-conrp} (2) and Theorem \ref{thm-FinHopftri-CH}] Let $H$ be an affine Hopf algebra over an algebraically closed field $\mathbbm{k}$ of characteristic zero with a central Hopf subalgebra  $C$ such that $H$ is a finitely generated $C$-module. Then the following hold.\label{intro-thmA}
\begin{itemize}
\item[(1)] As a $C$-module, $H$ is a finitely generated projective module of constant rank, say $d$.
\item[(2)] $(H,C,\text{tr}_{HS})$ is a Cayley-Hamilton Hopf algebra of degree $d$, where $\text{tr}_{HS}: H\to C$ is the Hattori-Stallings trace map of $H$ viewed as a projective $C$-module. 
\end{itemize}
\end{theorem}

See \S \ref{trHS-CHHOPF-sec} for the definition of the Hattori-Stallings trace map. Indeed,  Theorem \ref{intro-thmA} provides numerous examples of Cayley-Hamilton Hopf algebras, including many quantum groups at roots of unity and many AS-regular Hopf algebras of low GK-dimension.


In the rest of the introduction, we assume that the pair $(H,C)$ of Hopf algebras satisfies the conditions of Theorem \ref{intro-thmA}.
Note that $\text{maxSpec}C$ is an affine algebraic group with respect to convolution (see \S \ref{trHS-CHHOPF-sec}). Thus, we obtain \textit{a family of finite-dimensional $\mathbbm{k}$-algebras parameterized by the affine algebraic group} $\text{maxSpec}C$:
$$
\mathcal{B}_{H}\coloneqq\{H/\mathfrak{m}H\mid \mathfrak{m}\in \text{maxSpec}C\}.
$$
Denote the counit of $H$ by $\varepsilon$. Then $\overline{\varepsilon}\coloneqq \varepsilon|_{C}$ is the counit of the central Hopf subalgebra $C$, and $\mathfrak{m}_{\overline{\varepsilon}}\coloneqq \text{Ker}\overline{\varepsilon}\in \text{maxSpec}C$ is the augmentation ideal of $C$, which is also the identity element of the algebraic group $\text{maxSpec}C$. The coproduct of $H$ naturally endows the \textit{identity fiber algebra} $H/\mathfrak{m}_{\overline{\varepsilon}}H$ with a finite-dimensional Hopf algebra structure, and endows each fiber algebra $H/\mathfrak{m}H$ with both a left and a right $H/\mathfrak{m}_{\overline{\varepsilon}}H$-comodule algebra structure, see \S \ref{trHS-CHHOPF-sec}. Thus, the identity fiber algebra $H/\mathfrak{m}_{\overline{\varepsilon}}H$ can be regarded as a distinguished member in $\mathcal{B}_{H}$.


In the major part of this paper, we investigate the impact of the identity fiber algebras on the representation theory of such a class of Hopf algebras, especially how the \textit{Chevalley property} (which means the tensor product of any
two finite-dimensional irreducible modules is completely reducible) of the identity fiber algebra relates to the representation theory and \textit{discriminant ideals} of the Hopf algebra in question by applying the theories of tensor categories and Cayley-Hamilton algebras.

It was proved in \cite[Thm. 2.9 (8)]{MR4201485} that the global dimension of $H$ is finite if the identity fiber algebra $H/\mathfrak{m}_{\overline{\varepsilon}}H$ is semisimple. Inspired by this result, taking a representation-theoretic point of view, we naturally arrive at the following problem:

\medskip
\noindent
{\bf{Problem.}} For an affine Hopf algebra $H$ with a central Hopf subalgebra $C$ such that $H$ is a finitely generated $C$-module, determine how the structure of the identity fiber algebra $H/\mathfrak{m}_{\overline{\varepsilon}}H$ affects the representation theory of $H$.
\medskip

As mentioned before, the irreducible representations of $H$ can be studied through the family of fiber algebras $\mathcal{B}_{H}$, and each fiber algebra $H/\mathfrak{m}H$ is both a left and a right comodule algebra over the identity fiber Hopf algebra $H/\mathfrak{m}_{\overline{\varepsilon}}H$. Let $H/\mathfrak{m}H\text{-mod}$ be the category of finite-dimensional modules over $H/\mathfrak{m}H$ for any $\mathfrak{m}\in \text{maxSpec}C$. Then $H/\mathfrak{m}_{\overline{\varepsilon}}H\text{-mod}$ is a finite tensor category, and for each fiber algebra $H/\mathfrak{m}H$, $H/\mathfrak{m}H\text{-mod}$ is both a left and a right module category over the tensor category $H/\mathfrak{m}_{\overline{\varepsilon}}H\text{-mod}$, see \S \ref{trHS-CHHOPF-sec}. Therefore, studying the relationship between $H/\mathfrak{m}H\text{-mod}$ and $H/\mathfrak{m}_{\overline{\varepsilon}}H\text{-mod}$ is essential for addressing the aforementioned problem. This leads us to prove that \textit{for any} $\mathfrak{m}\in \text{maxSpec}C$\textit{, as a module category over the finite tensor category }$H/\mathfrak{m}_{\overline{\varepsilon}}H\text{-mod}$\textit{,} $H/\mathfrak{m}H\text{-mod}$ \textit{is an indecomposable exact module category} (Proposition \ref{thm-indexacmocat-fiberalg}). This result has two direct consequences.
\begin{itemize}
\item[(1)] 
 For any $\mathfrak{m}\in \text{maxSpec}C$, the Grothendieck group $\text{Gr}(H/\mathfrak{m}H)$ of $H/\mathfrak{m}H\text{-mod}$ can be naturally viewed as a $\mathbb{Z}_{+}$-module over the Grothendieck ring $\text{Gr}(H/\mathfrak{m}_{\overline{\varepsilon}}H)$ of the finite tensor category $H/\mathfrak{m}_{\overline{\varepsilon}}H\text{-mod}$ under the tensor product action (see  \S \ref{sec-H/mHmocat-exact}). According to \cite[Prop. 7.7.2]{MR3242743}, the fact that $H/\mathfrak{m}H\text{-mod}$ is an indecomposable exact module category over $H/\mathfrak{m}_{\overline{\varepsilon}}H\text{-mod}$ implies that $\text{Gr}(H/\mathfrak{m}H)$ is an irreducible $\mathbb{Z}_{+}$-module over $\text{Gr}(H/\mathfrak{m}_{\overline{\varepsilon}}H)$. The proofs of Theorems \ref{intro-thmB} and \ref{intro-thmC} are based on this observation, which also allows us to give new treatments to some old results. For example, when the identity fiber algebra $H/\mathfrak{m}_{\overline{\varepsilon}}H$ is basic, it was proved in \cite[Prop. III.4.11]{MR1898492} that the group action of the group of characters $G((H/\mathfrak{m}_{\overline{\varepsilon}}H)^{\circ})$ on the set $\text{Irr}(H/\mathfrak{m}H)$ of the isomorphism classes of irreducible $H/\mathfrak{m}H$-modules via tensor product is transitive for all $\mathfrak{m}\in \text{maxSpec}C$ by using Hopf algebra methods. We provide a \textit{tensor-categorical} proof of this result through the above observation (Remark \ref{rmk-irreZ+-appnewpf}). 

\item[(2)] The exactness of the module category $H/\mathfrak{m}H\text{-mod}$ over the tensor category $H/\mathfrak{m}_{\overline{\varepsilon}}H\text{-mod}$ enables us to prove that if $H$ is \textit{involutory} then $H$ has the Chevalley property (Corollary \ref{cor-invrHop-chev}). There has been extensive research on finite-dimensional Hopf algebras with the Chevalley property or dual Chevalley property, see \cite{MR1852304,MR2037722,MR4396654,MR4801656} for example. To our knowledge, there is little research in the infinite-dimensional case.
\end{itemize}

Another approach to explore the above problem is to consider the \textit{$\overline{\varepsilon}$-Chevalley locus} of the pair $(H,C)$,
defined by
\begin{align*}
\text{Chev}_{\overline{\varepsilon}}(H,C)\coloneqq \{\mathfrak{m}\in \text{maxSpec}C \mid &V\otimes W  \text{\ is a completely reducible\ } H/\mathfrak{m}H\text{-module for all\ }\\
&[V]\in \text{Irr}(H/\mathfrak{m}_{\overline{\varepsilon}}H)\text{ and\ }[W]\in \text{Irr}(H/\mathfrak{m}H)\}.
\end{align*}
By definition, $\text{Chev}_{\overline{\varepsilon}}(H,C)$ consists of points in $\text{maxSpec}C$ corresponding to the fiber algebras whose completely reducible representations remain completely reducible under the tensor action with the irreducible representations of the identity fiber algebra. If $H$ has the Chevalley property, then $\text{Chev}_{\overline{\varepsilon}}(H,C)=\text{maxSpec}C$. The converse is not true, see Example \ref{eg-cheloc=maxSp-notChprop}. 

Surprisingly,  $\text{Chev}_{\overline{\varepsilon}}(H,C)=\text{maxSpec}C$ if and only if the identity fiber algebra $H/\mathfrak{m}_{\overline{\varepsilon}}H$ has the Chevalley property as proved in Proposition \ref{prop-epiChev-mepiC}. 
Therefore, the $\overline{\varepsilon}$-Chevalley locus of the pair $(H,C)$ measures how far the identity fiber algebra $H/\mathfrak{m}_{\overline{\varepsilon}}H$ is from having the Chevalley property.

Brown-Yakimov \cite{MR3886192} discovered a deep relation between the dimensions of the irreducible representations of a Cayley-Hamilton algebra and its \textit{discriminant ideals} \cite{MR1972204} and \textit{modified discriminant ideals} \cite{MR3415697}. Discriminants (or more generally discriminant ideals) are powerful tools in both commutative and non-commutative algebras. In recent years, discriminants of non-commutative algebras and their applications have been intensely studied, see for example \cite{MR3281142,MR3415697,MR3644228,MR3663593,MR3720799,MR3614142,MR4413274}. However, little is currently known about the discriminant ideals. For an algebra with trace $(A,C,\text{tr})$ and a positive integer $k$, its $k$-discriminant ideal and modified $k$-discriminant ideal (Definition \ref{def-discriid}) are denoted by
$$D_{k}(A/C;\text{tr}), \text{ and\ } MD_{k}(A/C;\text{tr}),$$
respectively. These are ideals of the central subalgebra $C$ of $A$, satisfying
\begin{equation}\label{eq-mdcond}
MD_{k}(A/C;\text{tr})\supseteq D_{k}(A/C;\text{tr}),\ \text{for any positive integer\ }k.
\end{equation}
 Moreover, the modified discriminant ideals give a descending chain of ideals of $C$:
\begin{equation}\label{eq-desidch-MD}
MD_{1}(A/C;\text{tr})\supseteq MD_{2}(A/C;\text{tr}) \supseteq \cdots\supseteq MD_{k}(A/C;\text{tr})\supseteq MD_{k+1}(A/C;\text{tr})\supseteq \cdots.
\end{equation}
Let $(A,C,\text{tr})$ be an affine Cayley-Hamilton algebra over an algebraically closed field $\mathbbm{k}$ of characteristic zero. By \cite[Thm. 4.5 (a)]{MR1288995}, $A$ is a finitely generated $C$-module. Consequently, the chain \eqref{eq-desidch-MD} has only finitely many nonzero terms, that is, there is a positive integer $N$ such that
\begin{equation}\label{eq-chain-fini0tem}
0=MD_{N+1}(A/C;\text{tr})=MD_{N+2}(A/C;\text{tr})=MD_{N+3}(A/C;\text{tr})=\cdots.
\end{equation}
By \cite[Thm. 4.1 (b)]{MR3886192}, for any positive integer $k$,
\begin{equation}\label{eq-discczero-BY}
\begin{aligned}
\mathcal{V}_{k}\coloneqq \mathcal{V}(D_{k}(A/C;\text{tr}))&=\mathcal{V}(MD_{k}(A/C;\text{tr}))\\
&=\Big\{\mathfrak{m}\in \text{maxSpec}C\mid \sum_{[V]\in \text{Irr}(A/\mathfrak{m}A)}(\dim_{\mathbbm{k}}V)^{2} < k \Big\}.
\end{aligned} 
\end{equation}
Let $W_{l}(G(H^{\circ}))$ and $ W_{r}(G(H^{\circ}))$ be the left and right winding automorphism groups of $H$, respectively. See \S \ref{subsec-orbiwinda} for the background on winding automorphisms of Hopf algebras. It follows from \eqref{eq-discczero-BY} that for any affine Cayley-Hamilton Hopf algebra $(H,C,\text{tr})$, each zero locus of a discriminant ideal of $(H,C,\text{tr})$ is a disjoint union of some orbits of $\text{maxSpec}C$ under the action of $W_{l}(G(H^{\circ}))$ or $W_{r}(G(H^{\circ}))$. Therefore, studying the orbits of $\text{maxSpec}C$ under the action of the left (or right) winding automorphism group is fundamental for understanding the structure of the zero loci of discriminant ideals of $(H,C,\text{tr})$ and the irreducible representations of the fiber algebras of $H$.

This motivates our second result, which says that for any affine Cayley-Hamilton Hopf algebra such that its identity fiber algebra has the Chevalley property, all non-empty zero loci of its discriminant ideals contain a common orbit of the winding automorphism group action.



\begin{theorem}
[Theorem \ref{thm-Chein-VDneth-cwao}] Let $(H,C,\text{tr})$\label{intro-thmB} be an affine Cayley-Hamilton Hopf algebra over an algebraically closed field $\mathbbm{k}$ of characteristic zero such that the identity fiber algebra $H/\mathfrak{m}_{\overline{\varepsilon}}H$ has the Chevalley property. Then for any non-empty zero locus $\mathcal{V} $ of a discriminant ideal of $(H,C,\text{tr})$,
$$\mathcal{V} \supseteq W_{l}(G(H^{\circ}))(\mathfrak{m}_{\overline{\varepsilon}})= W_{r}(G(H^{\circ}))(\mathfrak{m}_{\overline{\varepsilon}}).$$
\end{theorem}

The proof of Theorem \ref{intro-thmB}  is based on two aforementioned results: first, for each $\mathfrak{m}\in \text{maxSpec}C$, $\text{Gr}(H/\mathfrak{m}H)$ is an irreducible $\mathbb{Z}_{+}$-module over $\text{Gr}(H/\mathfrak{m}_{\overline{\varepsilon}}H)$; secondly, $\text{Chev}_{\overline{\varepsilon}}(H,C)=\text{maxSpec}C$ under the assumptions of Theorem \ref{intro-thmB}. These facts allow us to apply Frobenius-Perron theory to study the behavior of the square dimension function (Theorem \ref{thm-sdchevpo-geqmepi}), thereby obtaining the result of Theorem \ref{intro-thmB}. By the proof of \cite[Thm. C (a)]{mi2025lowest}, the orbit
$$\mathcal{I}\coloneqq W_{l}(G(H^{\circ}))(\mathfrak{m}_{\overline{\varepsilon}})= W_{r}(G(H^{\circ}))(\mathfrak{m}_{\overline{\varepsilon}})$$ is a subgroup of $\text{maxSpec}C$. It is proved in Theorem \ref{thm-leftcossta-isofiber-or} that for any $\mathfrak{m}\in \text{maxSpec}C$, the left (right) coset of $\mathcal{I}$ in $\text{maxSpec}C$ containing $\mathfrak{m}$ is precisely the orbit of $\mathfrak{m}$ of the right (respectively, left) winding automorphism group action, see also Remark \ref{rmk-ricoscas-wnorsub1}.


In the special case of affine Cayley-Hamilton Hopf algebras where the identity fiber algebras are basic, the conclusion of Theorem \ref{intro-thmB} can be derived from \cite[Thm. C (b)]{mi2025lowest}. There are many important affine Hopf algebras admitting large central Hopf subalgebras such that the corresponding identity fiber algebras are basic. For example, the quantized coordinate rings at roots of unity and the big quantized Borel subalgebras at roots of unity \cite[Thms. 5.5 and 5.7]{mi2025lowest}. Nevertheless, the approach in \cite{mi2025lowest} relies heavily on the assumption that the identity fiber algebra is basic. Moreover, even for the group algebras of central extensions of finite groups by finitely generated abelian groups, the corresponding identity fiber algebras have the Chevalley property but are not basic in general (Example \ref{eg-grpalg-cenex1}).


In \S \ref{sec-app}, we present two applications. One concerns the \textit{lowest discriminant ideal} of Cayley-Hamilton Hopf algebras (Theorem \ref{intro-thmC}), the other is related to the Chevalley property of Cayley-Hamilton Hopf algebras (Theorem \ref{intro-thmD}). 

A discriminant ideal or a modified discriminant ideal of  an algebra with trace $(A,C,\text{tr})$ is said to be \textit{trivial} if it is $0$ or $C$. For an affine Cayley-Hamilton algebra $(A,C,\text{tr})$, by \eqref{eq-mdcond} and \eqref{eq-chain-fini0tem}, there are only finitely many non-trivial (modified) discriminant ideals. Thus, there is a unique positive integer $\ell$ for which
$$\varnothing=\mathcal{V}_{1}=\cdots=\mathcal{V}_{\ell-1}\subsetneq \mathcal{V}_{\ell}.$$
The corresponding discriminant ideal $D_{\ell}(A/C;\text{tr})$ (modified discriminant ideal $MD_{\ell}(A/C;\text{tr})$) is called the \textit{lowest  discriminant ideal} (respectively, \textit{the lowest modified discriminant ideal}) \cite{mi2025lowest}. In \cite[Thm. B]{mi2025lowest}, the authors described the zero loci of the lowest discriminant ideals of affine Cayley-Hamilton Hopf algebras $(H,C,\text{tr})$ whose identity fiber algebras are basic, and showed that in this setting, the lowest discriminant ideal of $(H,C,\text{tr})$ is of level
$$|G((H/\mathfrak{m}_{\overline{\varepsilon}}H)^{\circ})|+1,$$
where $G((H/\mathfrak{m}_{\overline{\varepsilon}}H)^{\circ})$ is the group of characters of the finite-dimensional Hopf algebra $H/\mathfrak{m}_{\overline{\varepsilon}}H$. 
The first application is to generalize some of the related results to a more general situation.

\begin{theorem}
[Theorem \ref{thm-lowCHHopf-minlif1}] Let $(H,C,\text{tr})$\label{intro-thmC} be an affine Cayley-Hamilton Hopf algebra over an algebraically closed field $\mathbbm{k}$ of characteristic zero such that the identity fiber algebra $H/\mathfrak{m}_{\overline{\varepsilon}}H$ has the Chevalley property. Let $\text{Irr}(H/\mathfrak{m}_{\overline{\varepsilon}}H)=\{[V_{1}],[V_{2}],...,[V_{m}]\}$. Then the following hold.
\begin{itemize}
\item[(1)] The lowest discriminant ideal of $(H,C,\text{tr})$ is of level 
$$\ell = \text{FPdim}(\text{Gr}(H/\mathfrak{m}_{\overline{\varepsilon}}H))+1,$$
where $\text{FPdim}(\text{Gr}(H/\mathfrak{m}_{\overline{\varepsilon}}H))$ is the Frobenius-Perron dimension of $\text{Gr}(H/\mathfrak{m}_{\overline{\varepsilon}}H)$.
\item[(2)] Suppose $\mathfrak{m}\in \text{maxSpec}C$. The following are equivalent.
\begin{itemize}
\item[(i)] $\mathfrak{m}\in \mathcal{V}_{\ell}$; 
\item[(ii)] For any irreducible $H/\mathfrak{m}H$-module $W$, as $H$-modules
$$W\otimes W^{*}\cong \bigoplus_{i=1}^{m}V_{i}^{\oplus \dim_{\mathbbm{k}}\text{Hom}_{H}(V_{i}\otimes W,W)};$$
\item[(iii)] For any irreducible $H/\mathfrak{m}H$-module $W$, $W\otimes W^{*}$ is a completely reducible $H$-module.
\end{itemize}
\end{itemize}
\end{theorem}

Note that if $H/\mathfrak{m}_{\overline{\varepsilon}}H$ is basic, then $\text{Gr}(H/\mathfrak{m}_{\overline{\varepsilon}}H)\cong \mathbb{Z}[G((H/\mathfrak{m}_{\overline{\varepsilon}}H)^{\circ})]$  as $\mathbb{Z}_{+}$-rings (Remark \ref{rmk-irreZ+-appnewpf}). It follows that $\text{FPdim}(\text{Gr}(H/\mathfrak{m}_{\overline{\varepsilon}}H))=\text{FPdim}(\mathbb{Z}[G((H/\mathfrak{m}_{\overline{\varepsilon}}H)^{\circ})])=|G((H/\mathfrak{m}_{\overline{\varepsilon}}H)^{\circ})|$.


The second application provides a necessary condition for a Cayley-Hamilton Hopf algebra  having the Chevalley property by using the discriminant ideals.


\begin{theorem}
[Theorem \ref{thm-CHHopf-Chev-dis}] Let $(H,C,\text{tr})$ be an affine Cayley-Hamilton Hopf algebra over an algebraically closed field $\mathbbm{k}$ of characteristic zero. If $H$ has the Chevalley property, then all the discriminant ideals of $(H,C,\text{tr})$ are trivial.\label{intro-thmD}
\end{theorem}

Thus, for any affine Hopf algebra which is a finitely generated free module over some central Hopf subalgebra, the discriminant can serve as a test for whether the Hopf algebra in question has the Chevalley property. We are able to provide some necessary and sufficient conditions for affine prime Cayley-Hamilton Hopf algebras to have the Chevalley property (Corollary \ref{cor-chev-primiff}). In fact, any affine prime Cayley-Hamilton Hopf algebra has the Chevalley property if and only if it is commutative. It follows that many Cayley-Hamilton Hopf algebras in the setting of Theorem \ref{intro-thmC} do not possess the Chevalley property. This includes many important families of quantum groups at roots of unity.

\begin{example}[the quantized coordinate ring $\mathcal{O}_{\epsilon}(G)$]
Let $G$ be a simply connected complex semisimple algebraic group of positive dimension with  Lie algebra $\mathfrak{g}$, let $\ell\geq 3$ be an odd positive integer, coprime to $3$ if $\mathfrak{g}$ contains a factor of type $G_{2}$, and let $\epsilon\in \mathbb{C}$ be a primitive $\ell$-th root of unity.  Consider the quantized coordinate ring $\mathcal{O}_{\epsilon}(G)$ over $\mathbb{C}$ \cite{MR1296515}. By \cite[Thm. 7.2]{MR1296515} and Theorem \ref{intro-thmA}, there is a central Hopf subalgebra $C_{\epsilon}(G)$ such that $(\mathcal{O}_{\epsilon}(G),C_{\epsilon}(G),\text{tr}_{HS})$ is an affine Cayley-Hamilton Hopf algebra of degree $\ell^{\dim G}$. Moreover, by \cite[Cor. 7.3 and Thm. 7.4]{MR1296515} and \cite[Prop. (2) III. 7.7]{MR1898492}, $\mathcal{O}_{\epsilon}(G)$ is a domain that is not commutative, and the corresponding identity fiber algebra $\mathcal{O}_{\epsilon}(G)/\mathfrak{m}_{\overline{\varepsilon}}\mathcal{O}_{\epsilon}(G)$ is basic. This shows that 
the quantized coordinate ring $\mathcal{O}_{\epsilon}(G)$ does not have the Chevalley property even though it is in the setting of Theorem \ref{intro-thmC}.
\end{example}
 \begin{remark}
In contrast, when the quantum parameter $q\in \mathbb{C}^{\times}$ is not a root of unity, the quantized enveloping algebras $\mathcal{U}_{q}(\mathfrak{g})$, the quantized Borel subalgebras $\mathcal{U}_{q}^{\geq 0}(\mathfrak{g})$, and the quantized coordinate rings $\mathcal{O}_{q}(G)$, all defined over $\mathbb{C}$, have the Chevalley property. The case of quantized enveloping algebras is well-known, while those of the latter two follow directly from \cite[Prop. 4.4 and Thm. 11.5]{MR1297159}.
 \end{remark}

The paper is organized as follows. In \S \ref{sec-prelimin}, we briefly review some preliminaries regarding module-finite algebras, Cayley-Hamilton algebras, $\mathbb{Z}_{+}$-rings and the Frobenius-Perron dimensions. In \S \ref{sec-fiberalgsbasip}, we study some basic properties of the family of fiber algebras of affine Hopf algebras that admit large central Hopf subalgebras, and prove Theorem \ref{intro-thmA} and Proposition \ref{thm-indexacmocat-fiberalg}. In \S \ref{sec-Chevloc-disccon}, we first prove Proposition \ref{prop-epiChev-mepiC} and study the behavior of the square dimension function of a module-finite Hopf algebra (Theorem \ref{thm-sdchevpo-geqmepi}), then we prove Theorem \ref{intro-thmB}. In \S \ref{sec-app}, we prove Theorems \ref{intro-thmC} and \ref{intro-thmD} as applications. 

\section{Preliminaries}\label{sec-prelimin}
Throughout the paper, unless otherwise specified, $\mathbbm{k}$ is a fixed base field; all vector spaces are over $\mathbbm{k}$, and the term ``algebra'' always refers to an associative $\mathbbm{k}$-algebra with identity. An unadorned $\otimes$ signifies $\otimes_{\mathbbm{k}}$.  All modules considered here are left modules. For an algebra $A$, the category of all finite-dimensional $A$-modules is denoted by $A\text{-mod}$. The set of isomorphism classes of irreducible $A$-modules is denoted by $\text{Irr}(A)$. The center of $A$ is denoted by $Z(A)$, and the maximal spectrum of $A$ is denoted by $\text{maxSpec}A$. For any $A$-module $W$ of finite length and any irreducible $A$-module $V$, $[W:V]$ is the multiplicity of $V$ in a composition series of $W$. 


The coproduct, counit, and antipode of a Hopf algebra are usually denoted by $\Delta,\varepsilon$ and $S$, respectively. For any Hopf algebra $H$ and $h\in H$, we use Sweedler's notation $\Delta(h)=\sum h_{(1)}\otimes h_{(2)} \in H \otimes H$. Let $H^{\circ}$ be the finite dual Hopf algebra of $H$. The set of group-like elements of $H$ is denoted by $G(H)$. Clearly, $G(H^{\circ})$ coincides with the set of all characters of $H$. Canonically, there is a  one-to-one correspondence between $G(H^{\circ})$ and the isomorphism classes of $1$-dimensional $H$-modules. We refer the reader to \cite{MR1243637} for basic material on Hopf algebras.

In this section, we provide some preliminaries on module-finite algebras, Cayley-Hamilton algebras, $\mathbb{Z}_{+}$-rings and the Frobenius-Perron dimensions. For more detailed background, we refer the reader to our primary references \cite{MR1898492,MR1811901,MR883869,MR1288995,MR2178656,MR3242743}.


\subsection{Module-finite algebras}\label{sec-mofialgs}
Let $A$ be an algebra with a central subalgebra $C$ such that $A$ is a finitely generated $C$-module. Such an algebra $A$ is referred to as a \textit{module-finite $C$-algebra} \cite[p.xv]{MR2080008}. In this subsection, we always assume that $A$ is an affine module-finite $C$-algebra over  an algebraically closed field $\mathbbm{k}$. By the Artin-Tate Lemma, $C$ is an affine $\mathbbm{k}$-algebra (in particular, $C$ is Noetherian). Thus, for any $\mathfrak{m}\in \text{maxSpec}C$, $C/\mathfrak{m}\cong \mathbbm{k}$ by Hilbert's Nullstellensatz. There is a family of finite-dimensional $\mathbbm{k}$-algebras parameterized by the points in $\text{maxSpec}C$:
$$
\mathcal{B}_{A}\coloneqq\{A/\mathfrak{m}A\mid \mathfrak{m}\in \text{maxSpec}C\}.
$$
If $A$ can be generated by $d$ elements as a $C$-module, then each member in $\mathcal{B}_{A}$ is a finite-dimensional $\mathbbm{k}$-algebra with $\mathbbm{k}$-dimension $\leqslant d$. If $A$ is further a finitely generated projective $C$-module of constant rank $d$ (for example, $A$ is a free $C$-module of rank $d$), then by Nakayama's Lemma, all algebras in $\mathcal{B}_{A}$ have $\mathbbm{k}$-dimension $d$.

Note that $A$ is an affine PI algebra over $\mathbbm{k}$.  Any irreducible $A$-module is finite-dimensional over $\mathbbm{k}$ \cite[Thm. 13.10.3 (i)]{MR1811901}. It follows from  Kaplansky's Theorem that there is a bijection between the set of isomorphism classes  $\text{Irr}(A)$ of irreducible $A$-modules and the set of maximal ideals $\text{maxSpec}A $ of $A$ given by $[V]\mapsto \text{Ann}_{A}V$. 

By \cite[Thm. 13.8.14]{MR1811901}, 
$\pi:\text{maxSpec}A\to \text{maxSpec}C, Q\mapsto Q\cap C$ is a  surjective map. 
For any irreducible $A$-module $V$, $V$ is an irreducible module over  $A/\mathfrak{m}A$, where $\mathfrak{m}=(\text{Ann}_{A}V)\cap C \in \text{maxSpec}C$. 
Therefore, studying the irreducible representations of $A$ is equivalent to studying the irreducible representations of all algebras in $\mathcal{B}_{A}$ in some sense. For any $\mathfrak{m}\in \text{maxSpec}C$, the set of irreducible $A$-modules corresponding to the fiber $\pi^{-1}(\mathfrak{m})$ of $\mathfrak{m}$ under $\pi$ is precisely $\text{Irr}(A/\mathfrak{m}A)$. The finite-dimensional algebra $A/\mathfrak{m}A$ is called the \textit{fiber algebra} at $\mathfrak{m}$.

\subsection{Cayley-Hamilton algebras}\label{sec-CH algs}
Throughout this subsection, $\mathbbm{k}$ is of characteristic zero. Suppose  $\text{M}_{n}(C)$ is the $n\times n$ matrix algebra over a commutative $\mathbbm{k}$-algebra $C$. 
Let
$$p_{n,X}(t)=t^{n}-c_{1}(X)t^{n-1}+\cdots+(-1)^{n-1}c_{n-1}(X)t+(-1)^{n}c_{n}(X)\in C[t]$$ be the 
characteristic polynomial of $X\in \text{M}_{n}(C)$. It is well-known (see for example \cite[Cor. 3.3.3 (i)]{MR4249615}) that the coefficient $c_{k}(X)$ ($1\leq k\leq n$)  is of the following form
in terms of the traces of powers of $X$:
\begin{equation}\label{eq-trfor-chapolym}
c_{k}(X)=\frac{1}{k!}\left|\begin{array}{ccccc}\text{tr}(X) & 1 & 0 & \cdots & 0 \\ \text{tr}(X^{2}) & \text{tr}(X) & 2 & \cdots & 0 \\ \vdots & \vdots & \vdots & & \vdots \\ \text{tr}(X^{k-1}) & \text{tr}(X^{k-2}) & \text{tr}(X^{k-3}) & \cdots & k-1 \\ \text{tr}(X^{k}) & \text{tr}(X^{k-1}) & \text{tr}(X^{k-2}) & \cdots & \text{tr}(X)\end{array}\right|,
\end{equation}
where $\text{tr}(X^{i})$ is the usual trace of the matrix $X^{i}$. 

Now, for any algebra with trace $(A,C,\text{tr})$ over $\mathbbm{k}$, any positive integer $k$, and any element $a\in A$, one can formally define $c_{k}(a)\in C$ via formula \eqref{eq-trfor-chapolym}. Then, for any positive integer $n$, 
\begin{equation}\label{eq-charpoly-CHalg}
p_{n,a}(t)=t^{n}-c_{1}(a)t^{n-1}+\cdots+(-1)^{n-1}c_{n-1}(a)t+(-1)^{n}c_{n}(a)\in C[t].
\end{equation}
is called the \textit{$n$-th characteristic polynomial} of $a \in A$. 

\begin{definition}
[\cite{MR883869,MR1288995}]Let $(A,C,\text{tr})$ be an algebra with trace over a field $\mathbbm{k}$ of characteristic zero and $n$ a positive integer. The algebra with trace $(A,C,\text{tr})$ is called a Cayley-Hamilton algebra of degree $n$ if $p_{n,a}(a)=0$ for all $a\in A$ and $\text{tr}(1)=n$.
\end{definition}

Clearly, if $\text{tr}:\text{M}_{n}(C)\to C$ is the usual trace map of the matrix algebra, then \eqref{eq-charpoly-CHalg} coincides with the usual characteristic polynomial of a matrix, and  $(\text{M}_{n}(C),C,\text{tr})$ is a Cayley-Hamilton algebra of degree $n$.
\begin{example}
[regular trace]Let \label{eg-regulartr} $A$ be a module-finite $C$-algebra such that $A$ is a free $C$-module  of rank $n$. Suppose $\{b_{1},...,b_{n}\}$ is a $C$-basis of $A$. Then for any element $a\in A$, there is a unique matrix $\iota(a)\coloneqq{(c_{ij})}\in \text{M}_{n}(C)$ such that $a b_j=\sum_i c_{ij} b_i$ for all $ 1\le j\le n$. This defines an algebra embedding $\iota:A\to\text{M}_{n}(C)$. Then the composition of $\iota$ with the usual trace map on $\text{M}_{n}(C)$ gives a trace map $\text{tr}_{\text{reg}}:A\to C$. That is, we have the following commutative diagram:
$$\begin{tikzcd}
A \arrow[r, "\iota"] \arrow[rd, "\text{tr}_{\text{reg}}"'] & \text{M}_{n}(C) \arrow[d, "\text{tr}"] \\
 & {C} 
\end{tikzcd}$$
where $\text{tr}:\text{M}_{n}(C)\to C$ denotes the usual trace map. Clearly, the trace map $\text{tr}_{\text{reg}}$ is independent of the choice of the $C$-basis of $A$, which is called the \textit{regular trace} of $A$. It is straightforward to verify that the algebra with trace $(A,C,\text{tr}_{\text{reg}})$ is a Cayley-Hamilton algebra of degree $n$.

\end{example}
\begin{example}
[reduced trace]Let $A$ be a module-finite $C$-algebra\label{eg-redtr} where $C=Z(A)$ is a normal domain such that $A$ is prime ring of PI-degree $n$. By Posner's Theorem, $A\otimes_{C}\text{Frac}(C)$ is a central simple algebra. Let $F\supseteq \text{Frac}(C)$ be a splitting field for the central simple algebra $A\otimes_{C}\text{Frac}(C)$. Then there is an $F$-algebra isomorphism $h:A\otimes_{C}F\to\text{M}_{n}(F)$. Let $\text{tr}:\text{M}_{n}(F)\to F$ denote the usual trace on the matrix algebra. By \cite[p.113]{MR1972204}, $\text{tr}(h(A\otimes_{C}F))\subseteq \text{Frac}(C)$ and the trace map $\widehat{\text{tr}}_{\text{red}}: A\otimes_{C}\text{Frac}(C)\to \text{Frac}(C),x\mapsto \text{tr}(h(x)) $ is independent of the choices of $F$ and $h$. Since $C$ is normal, by \cite[Prop. 8.11 (b)(ii)]{MR1811901}, the restriction of $\widehat{\text{tr}}_{\text{red}}$ to $A$ gives a trace map $A \to C$. The resulting trace map is denoted by $\text{tr}_{\text{red}}: A \to C$, and it is called the \textit{reduced trace} of $A$. It is not hard to prove that the algebra with trace $(A,C,\text{tr}_{\text{red}})$ is a Cayley-Hamilton algebra of degree $n$.
\end{example}
Many important examples of module-finite algebras fall into the situation of Example \ref{eg-regulartr} or Example \ref{eg-redtr}. Let $A$ be a prime affine module-finite twisted Calabi-Yau algebra $A$ (see \cite[Def. 1.2 (ii)]{MR4413364} for the definition) of PI-degree $n$. It follows from \cite[Prop. 2.9]{MR2492474} that $A$ is homologically homogeneous in the sense of Brown-Hajarnavis \cite{MR719665}. Then the center $Z(A)$ of $A$ is normal by \cite[Thm. 6.1]{MR719665}. It follows that $(A,Z(A),\text{tr}_{\text{red}})$ is a Cayley-Hamilton algebra of degree $n$.

\begin{definition}
[\cite{MR2178656}] A Cayley-Hamilton Hopf algebra is a Cayley-Hamilton algebra $(H,C,\text{tr})$ such that  $H$ is a Hopf algebra and $C$ is a central Hopf subalgebra of $H$.\label{def-CHHopf}
\end{definition}

\subsection{\texorpdfstring{$\mathbb{Z}_{+}$-rings}{Z+rings} and the Frobenius-Perron dimensions}\label{sec-Z+FPdim}
In this subsection, we recall some basic results on $\mathbb{Z}_{+}$-rings and the Frobenius-Perron dimensions.
The readers are referred to \cite{MR3242743} for details. Let $\mathbbm{k}$ be an algebraically closed field.

Let $R$ be a ring that is free as a $\mathbb{Z}$-module and has a $\mathbb{Z}$-basis $\{b_{i}\}_{i\in I}$. The pair $(R,\{b_{i}\}_{i\in I})$ is called a \textit{$\mathbb{Z}_{+}$-ring} if for any $i,j\in I$, the structure constants $c_{ij}^{k}$ in the product $b_{i}b_{j}=\sum_{k\in I}c_{ij}^{k}b_{k}$ are all non-negative integers. In this case, $\{b_{i}\}_{i\in I}$ is referred to as the \textit{$\mathbb{Z}_{+}$-basis} of $R$. When no confusion arises, the $\mathbb{Z}_{+}$-ring $(R,\{b_{i}\}_{i\in I})$ is simply denoted by $R$. If the identity $1$ of a $\mathbb{Z}_{+}$-ring $R$ lies in its $\mathbb{Z}_{+}$-basis, then $R$ is called a \textit{unital $\mathbb{Z}_{+}$-ring}.

Given a finite $\mathbbm{k}$-linear abelian category $\mathcal{C}$, $\text{Irr}(\mathcal{C})$ denotes the set of isomorphism classes of irreducible/simple objects in $\mathcal{C}$, and $\text{Gr}(\mathcal{C})$ is the Grothendieck group of $\mathcal{C}$. Given a $\mathbbm{k}$-algebra $A$, we abbreviate $\text{Gr}(A\text{-mod})$ as $\text{Gr}(A)$. If $\mathcal{C}$ is further a finite tensor category, then $\text{Gr}(\mathcal{C})$ is a transitive unital $\mathbb{Z}_{+}$-ring with $\mathbb{Z}_{+}$-basis $\text{Irr}(\mathcal{C})$ \cite[Prop. 4.5.4]{MR3242743}. In particular, for any finite-dimensional Hopf algebra $H$, the Grothendieck ring $\text{Gr}(H)$ is a unital $\mathbb{Z}_{+}$-ring with $\text{Irr}(H)$ as its $\mathbb{Z}_{+}$-basis.

Let $(R,\{b_{i}\}_{i=1}^{m})$ be a transitive unital $\mathbb{Z}_{+}$-ring of finite rank, where $\{b_{i}\}_{i=1}^{m}$ is the $\mathbb{Z}_{+}$-basis of $R$. The Frobenius-Perron dimension of  $a\in R$ is denoted by $\text{FPdim}(a)$. For a finite-dimensional Hopf algebra $H$ and an irreducible $H$-module $V$, the Frobenius-Perron dimension of $[V]$ is $\dim_{\mathbbm{k}}V$. The following facts are needed later.
\begin{lemma}
Let $(R,\{b_{i}\}_{i=1}^{m})$ be a transitive unital $\mathbb{Z}_{+}$-ring.\label{lem-z+ring-regel}
\begin{itemize}
\item[(1)]\cite[Prop. 3.3.6 (2)]{MR3242743} There is a nonzero element $\mathscr{R}\in R_{\mathbb{C}}\coloneqq R\otimes_{\mathbb{Z}}\mathbb{C}$, unique up to multiplication by a complex number, such that for all $a\in R$ one has $a\mathscr{R} =\text{FPdim}(a)\mathscr{R}$. Such an element $\mathscr{R}$ is called the \textit{regular element} of $R$.
\item[(2)]\cite[Prop. 3.3.11]{MR3242743} If $R$ is further a fusion ring, then $\sum_{i=1}^{m}\text{FPdim}(b_{i})b_{i} $ is a regular element of $R$, which is called the \textit{canonical regular element} of $R$.
\end{itemize}
\end{lemma}
For a fusion ring $(R,\{b_{i}\}_{i=1}^{m})$, the \textit{Frobenius-Perron dimension} of $R$, denoted by $\text{FPdim}(R)$, is defined to be
$\text{FPdim}(\mathscr{R})=\sum_{i=1}^{m}\text{FPdim}(b_{i})^{2}$. If the Grothendieck ring $\text{Gr}(H)$ of a finite-dimensional Hopf algebra $H$ over an algebraically closed field is a fusion ring (e.g., when $H$ is semisimple, see also \cite[Prop. 4.9.1]{MR3242743}), then
\begin{equation}\label{eq-fus-FPdim}
\text{FPdim}(\text{Gr}(H))=\sum_{[V]\in \text{Irr}(H)}(\dim_{\mathbbm{k}}V)^{2}.
\end{equation}


\section{Fiber algebras of module-finite Hopf algebras}\label{sec-fiberalgsbasip}

In this section, we first study some basic properties of the fiber algebras of affine Hopf algebras $H$ that admit large central Hopf subalgebras $C$. We then proceed to prove Theorem \ref{intro-thmA} and Proposition \ref{thm-indexacmocat-fiberalg} in \S \ref{trHS-CHHOPF-sec} and \ref{sec-H/mHmocat-exact}. In \S \ref{subsec-orbiwinda}, we investigate the orbits of points in $\text{maxSpec}C$ under the action of the left/right winding automorphism group of $H$. Throughout this section, $\mathbbm{k}$ is an algebraically closed field. 

\subsection{Resource of Cayley-Hamilton Hopf algebras}\label{trHS-CHHOPF-sec}
In this subsection, we assume further that  $\mathbbm{k}$ is of characteristic zero.
For any affine Hopf algebra $H$ admitting a large central Hopf subalgebra $C$, we show that $H$ is a finitely generated projective $C$-module of constant rank (Proposition \ref{prop-equfibalg-conrp}). Hence the pair $(H,C)$ can be equipped with  the Hattori-Stallings trace map $\text{tr}_{HS}:H\to C$  (see \eqref{eq-def-HStrace} for definition) \cite{MR175950,MR202807}, so that $(H,C, \text{tr}_{HS})$  is a Cayley-Hamilton Hopf algebra (Theorem \ref{thm-FinHopftri-CH}). In particular, any PI Hopf triple in the sense of Brown \cite{MR1676211,MR1898492} is a Cayley-Hamilton Hopf algebra with respect to the Hattori-Stallings trace map. 

For the convenience, we make the following hypothesis:
\begin{align*}
\text{(FinHopf) }& H \textit{ is an affine Hopf algebra over  }\mathbbm{k}\textit{ with a } \textit{central Hopf subalgebra } C \\
& \textit{such that } H\textit{ is a finitely generated }C\textit{-module}.
\end{align*}

In this case, $C$ is an affine commutative Hopf algebra by the Artin-Tate Lemma. Since the field $\mathbbm{k}$ is of characteristic zero, $\text{maxSpec}C$ is a smooth affine variety
with trivial canonical bundle by Cartier's Theorem (see \cite[Thms. 11.4 and 11.6]{MR547117}) and  \cite[Cor. 11.3]{MR547117}. Hence $C$ is Calabi-Yau. Moreover, the affine Calabi-Yau variety $\text{maxSpec}C$ is an algebraic group with respect to the convolution product. There is a canonical bijection between the group of characters $G(C^{\circ})$ and $\text{maxSpec}C$ given by $\chi\mapsto \mathfrak{m}_{\chi}\coloneqq \text{Ker}\chi$. Under this correspondence, for any two points $\mathfrak{m}_{\chi}$ and $\mathfrak{m}_{\psi}$ in $\text{maxSpec}C$, where $\chi,\psi\in G(C^{\circ})$, the \textit{convolution product} of $\mathfrak{m}_{\chi}$ and $\mathfrak{m}_{\psi}$ is given by 
$$\mathfrak{m}_{\chi}\ast\mathfrak{m}_{\psi}\coloneqq \mathfrak{m}_{\chi\ast\psi},$$
where $\chi\ast\psi$ denotes the convolution product of the characters $\chi$ and $\psi$. For the counit $\varepsilon:H\to\mathbbm{k}$, $\overline{\varepsilon}\coloneqq \varepsilon|_{C}:C\to\mathbbm{k}$ is the counit of $C$. The identity of $\text{maxSpec}C$ with respect to the convolution product is $\mathfrak{m}_{\overline{\varepsilon}}$. The inverse of any $\mathfrak{m}_{\chi}\in \text{maxSpec}C$ with respect to the convolution product is given by $\mathfrak{m}_{\chi S}$, where $S$ is the antipode of $H$. For any $\mathfrak{m}\in \text{maxSpec}C$, the inverse of $\mathfrak{m}$ with respect to the convolution product is denoted by $\mathfrak{m}^{-1}$. It follows from the bijectivity of the antipode $S$ \cite[Cor. 2]{MR2271358} that $S(\mathfrak{m}^{-1})=\mathfrak{m}$. Hence the antipode $S$ of $H$ induces an algebra anti-isomorphism $\overline{S}:H/\mathfrak{m}^{-1}H \to H/\mathfrak{m}H$.

Therefore, the family of fiber algebras $\mathcal{B}_{H}=\{H/\mathfrak{m}H\mid \mathfrak{m}\in \text{maxSpec} C\}$ is parameterized by the affine algebraic group $\text{maxSpec}C$. We will see that all algebras in $\mathcal{B}_{H}$ have the same $\mathbbm{k}$-dimension in Proposition \ref{prop-equfibalg-conrp}. For any two characters $\chi,\psi\in G(C^{\circ})$, obviously  $(\chi\otimes \psi)(\mathfrak{m}_{\chi\ast\psi})=0$, and $\Delta(\mathfrak{m}_{\chi\ast\psi})\subseteq \mathfrak{m}_{\chi}\otimes C+C\otimes \mathfrak{m}_{\psi}$. Then the coproduct $\Delta:H\to H\otimes H$ induces an algebra homomorphism
\begin{equation}\label{eq-coproun-algfibers}
\overline{\Delta}:H/\mathfrak{m}_{\chi\ast\psi}H\to H/\mathfrak{m}_{\chi}H\otimes H/\mathfrak{m}_{\psi}H.
\end{equation}
Thus, the identity fiber algebra $H/\mathfrak{m}_{\overline{\varepsilon}}H$ becomes a finite-dimensional Hopf algebra. It is referred to as the \textit{restricted Hopf algebra} associated to the pair $(H,C)$ in \cite{MR1898492}. It follows from \eqref{eq-coproun-algfibers} also that each fiber algebra $H/\mathfrak{m}H$ is a left and right comodule algebra over $H/\mathfrak{m}_{\overline{\varepsilon}}H$:
\begin{equation}\label{eq-righlefcomalh9}
\overline{\Delta}:H/\mathfrak{m}H\to H/\mathfrak{m}_{\overline{\varepsilon}}H\otimes H/\mathfrak{m}H,\ \overline{\Delta}:H/\mathfrak{m} H\to H/\mathfrak{m}H \otimes H/\mathfrak{m}_{\overline{\varepsilon}}H.
\end{equation}
Thus, \eqref{eq-righlefcomalh9} turns $H/\mathfrak{m}H\text{-mod}$ into both a left and a right module category over the finite tensor category $H/\mathfrak{m}_{\overline{\varepsilon}}H\text{-mod}$. We shall return to this in \S \ref{sec-H/mHmocat-exact}.




Since $C$ is Noetherian and $H$ is finitely generated as $C$-module, it follows that $H$ is a Noetherian Hopf algebra. By \cite[Thm. 3.3]{MR1228767}, \textit{$H$ is projective over $C$}. In fact, under the hypothesis (FinHopf), it can be shown that $H$ is a Frobenius extension of $C$ by either a dualizing complex argument or a Hopf-Galois extension argument. We will not need this fact in this paper, the readers are referred to \cite[Cor. III.4.7]{MR1898492} and \cite[\S 2.5]{MR2379091} for details. 
 
Next we show that $ H$ is a projective $C$-module of constant rank.
\begin{proposition}
Let $(H,C)$ be a pair of Hopf algebras satisfying the hypothesis (FinHopf).\label{prop-equfibalg-conrp} 
\begin{itemize}
\item[(1)]\cite[Lem. 4.6]{MR4877094} For any $\mathfrak{m},\mathfrak{n}\in \text{maxSpec}C$, $\dim_{\mathbbm{k}}H/\mathfrak{m}H=\dim_{\mathbbm{k}}H/\mathfrak{n}H$.
\item[(2)]As a finitely generated projective $C$-module, $H$ has constant rank.
\end{itemize}
\end{proposition}
\begin{proof}
(1) 
We provide a simpler proof here. 
Let $\mu: H/\mathfrak{m}H\otimes H/\mathfrak{m}H \to H/\mathfrak{m}H$ be the induced map of the multiplication map of $H$. Let $\beta$ be  the composition of the following maps:
$$\begin{tikzcd}[column sep=small]
H/\mathfrak{m}H\otimes H/\mathfrak{m}H \arrow[r, "\text{id}\otimes\overline{\Delta}"] & H/\mathfrak{m}H\otimes H/\mathfrak{m}H\otimes H/\mathfrak{m}_{\overline{\varepsilon}}H \arrow[r, "\mu\otimes \text{id}"] & H/\mathfrak{m}H\otimes H/\mathfrak{m}_{\overline{\varepsilon}}H.
\end{tikzcd}$$
In fact, $\beta:H/\mathfrak{m}H\otimes H/\mathfrak{m}H\to H/\mathfrak{m}H\otimes H/\mathfrak{m}_{\overline{\varepsilon}}H$ is the map $\overline{h}\otimes \overline{g}\mapsto \sum_{(g)}\overline{hg_{(1)}}\otimes \overline{g_{(2)}}.$

Then, $\beta$ is invertible with the inverse  given by the following commutative diagram:
$$\begin{tikzcd}[column sep=large]
H/\mathfrak{m}H\otimes H/\mathfrak{m}_{\overline{\varepsilon}}H \arrow[d, "\text{id}_{H/\mathfrak{m}H}\otimes \overline{\Delta}"'] \arrow[rr, "\beta^{-1}", dashed] &  & H/\mathfrak{m}H\otimes H/\mathfrak{m}H                                                                             \\
H/\mathfrak{m}H\otimes H/\mathfrak{m}^{-1}H\otimes H/\mathfrak{m}H \arrow[rr, "\text{id}_{H/\mathfrak{m}H}\otimes \overline{S}\otimes \text{id}_{H/\mathfrak{m}H}"]             &  & H/\mathfrak{m}H\otimes H/\mathfrak{m}H\otimes H/\mathfrak{m}H \arrow[u, "\mu\otimes \text{id}_{H/\mathfrak{m}H}"']
\end{tikzcd}$$
That is, $\beta^{-1}: H/\mathfrak{m}H\otimes H/\mathfrak{m}_{\overline{\varepsilon}}H\to H/\mathfrak{m}H\otimes H/\mathfrak{m}H$ is given by $\overline{h}\otimes \overline{g}\mapsto \sum_{(g)}\overline{hS(g_{(1)})}\otimes \overline{g_{(2)}}.
$

 It follows that $(\dim_{\mathbbm{k}}H/\mathfrak{m}H)^{2}=(\dim_{\mathbbm{k}}H/\mathfrak{m}H)(\dim_{\mathbbm{k}}H/\mathfrak{m}_{\overline{\varepsilon}}H).$ Hence (1) holds.
\par(2) Let $\mathfrak{m}\in \text{maxSpec}\,C$. By Nakayama's Lemma, 
$$\text{rank}_{C_{\mathfrak{m}}}H_{\mathfrak{m}}=\dim_{C_{\mathfrak{m}}/\mathfrak{m}C_{\mathfrak{m}}}H_{\mathfrak{m}}/\mathfrak{m}H_{\mathfrak{m}}=\dim_{C/\mathfrak{m}}H/\mathfrak{m}H=\dim_{\mathbbm{k}}H/\mathfrak{m}H.$$
It follows from (1) that $\text{rank}_{C_{\mathfrak{m}}}H_{\mathfrak{m}}=\text{rank}_{C_{\mathfrak{n}}}H_{\mathfrak{n}}$ for all $\mathfrak{m},\mathfrak{n}\in \text{maxSpec}C$. Hence $_{C}H$ has constant rank.
\end{proof}
\begin{remark}
Under the hypothesis (FinHopf),  one can also show that for each fiber algebra $H/\mathfrak{m}H$,  
$\mathbbm{k}\subseteq H/\mathfrak{m}H$ is an $H/\mathfrak{m}_{\overline{\varepsilon}}H$-Galois extension, see \cite[Lem. 3.1]{MR4600057}. 
\end{remark}
\begin{remark}
There exist pairs $(H,C)$ satisfying the hypothesis (FinHopf) such that $H$ is not a free module over $C$, see \cite[Rmk. 2.4 (2)]{MR4201485} for example.
\end{remark}

\begin{example}
[\cite{MR1296515,MR1863398}]Let $G$ be\label{eg-quantizecoornri} a simply connected complex semisimple algebraic group of positive dimension with Lie algebra $\mathfrak{g}$, let $\ell\geq 3$ be an odd positive integer, coprime to $3$ if $\mathfrak{g}$ contains a factor of type $G_{2}$, and let $\epsilon\in \mathbb{C}$ be a primitive $\ell$-th root of unity. Consider the quantized coordinate ring $\mathcal{O}_{\epsilon}(G)$ over $\mathbb{C}$ defined in \cite{MR1296515}. In \cite[Prop. 6.4]{MR1296515}, De Concini-Lyubashenko constructed a central Hopf subalgebra $C_{\epsilon}(G)$ of $\mathcal{O}_{\epsilon}(G)$ and proved that $C_{\epsilon}(G)\cong \mathcal{O}(G)$ as Hopf algebras, where $\mathcal{O}(G)$ is the coordinate ring of $G$. Under the algebraic group isomorphism $\text{maxSpec}C_{\epsilon}(G)\cong G$, the point $\mathfrak{m}_{\overline{\varepsilon}}$ corresponds precisely to the identity element of $G$. De Concini-Lyubashenko also proved that $\mathcal{O}_{\epsilon}(G)$ is a projective module over $C_{\epsilon}(G)$ of rank $\ell^{\dim G}$ \cite[Thm. 7.2]{MR1296515}. In fact, it can be further shown that $\mathcal{O}_{\epsilon}(G)$ is a free module over the central Hopf subalgebra $C_{\epsilon}(G)$ \cite[Prop. 2.2]{MR1863398}. According to Example \ref{eg-regulartr},  $(\mathcal{O}_{\epsilon}(G),C_{\epsilon}(G),\text{tr}_{\text{reg}})$ is a Cayley-Hamilton Hopf algebra of degree $\ell^{\dim G}$. By \cite[Prop. (2) III. 7.7]{MR1898492} or \cite[Thm. 5.7 (a)]{mi2025lowest}, the identity fiber algebra $\mathcal{O}_{\epsilon}(G)/\mathfrak{m}_{\overline{\varepsilon}}\mathcal{O}_{\epsilon}(G)$ is basic. Let $N$ be the number of the positive roots of the complex semisimple Lie algebra $\mathfrak{g}$. By \cite[Cor. 7.3 and Thm. 7.4]{MR1296515}, $\mathcal{O}_{\epsilon}(G)$ is a domain of PI-degree $\ell^{N}$. It follows from Posner's Theorem that $\mathcal{O}_{\epsilon}(G)$ is not commutative. Moreover, by \cite[Thm. B (ii)]{MR1482982} and \cite[Thm. 0.2 (1)]{MR1938745}, $\mathcal{O}_{\epsilon}(G)$ is an affine AS-regular Hopf algebra.
\end{example}

Before we prove Theorem \ref{FinHopf-CH-Hopf} we need to introduce  Hattori-Stallings trace maps. 

Let $A$ be a module-finite $C$-algebra over $\mathbbm{k}$ such that $A$ is projective over $C$. Let $\{x_{1},...,x_{n}\}\subseteq A$ and $\{x_{1}^{*},...,x_{n}^{*}\}\subseteq \text{Hom}_{C}(A,C)$ be a dual basis of the projective module $_{C}A$. It is straightforward to verify that the map
\begin{equation}\label{eq-def-HStrace}
\text{tr}_{HS}:A\to C,a\mapsto \sum_{i=1}^{n}x_{i}^{*}(ax_{i})
\end{equation}
is a trace map, and it is independent of the choices of the dual basis $\{x_{1},...,x_{n}\}$ and $\{x_{1}^{*},...,x_{n}^{*}\}$ of $_{C}A$ (see \cite[\S 1.1]{MR4413274} for details). The trace map $\text{tr}_{HS}$  defined by \eqref{eq-def-HStrace} is called the \textit{Hattori-Stallings trace map} of $A$ over $C$ \cite{MR175950,MR202807}.  If $A$ is free over $C$, then $\text{tr}_{HS}$ coincides with the regular trace in Example \ref{eg-regulartr}. 

 Hence, under the the hypothesis (FinHopf), one can consider the Hattori-Stallings trace map $\text{tr}_{HS}:H\to C$. Moreover, by Proposition \ref{prop-equfibalg-conrp}(2), $_{C}H$ is a finitely generated projective module of constant rank. Now Theorem \ref{intro-thmA} follows from Proposition \ref{prop-equfibalg-conrp} (2) and the following.
 
\begin{theorem} \label{FinHopf-CH-Hopf}
Let $(H,C)$ be a pair of Hopf algebras satisfying the hypothesis (FinHopf). Then $(H,C,\text{tr}_{HS})$ is a Cayley-Hamilton Hopf algebra of degree $n$, where $n$ is the constant rank of $H$ as a projective $C$-module.\label{thm-FinHopftri-CH}
\end{theorem}
\begin{remark}
Let $(H,C)$ be a pair of Hopf algebras satisfying the hypothesis (FinHopf). If further $H$ is a prime ring, then the triple $(H,Z(H),C)$ is said to be a \textit{PI Hopf triple}. In this setting, the central Hopf subalgebra $C$, being a smooth domain (recall that $\text{char}\mathbbm{k}=0$), is of course normal. Then, one can apply \cite[Thm. 4.2]{MR2178656} to determine all the trace maps $\text{tr}:H\to C$ such that $(H,C,\text{tr})$ is a Cayley-Hamilton Hopf algebra.
\end{remark}

Theorem \ref{thm-FinHopftri-CH} is immediate from the following theorem.

\begin{theorem} 
Let $A$ be a module-finite $C$-algebra such that $A$ is a projective $C$-module of rank $n$. Then $(A,C,\text{tr}_{HS})$ is a Cayley-Hamilton algebra of degree $n$.
\end{theorem}
\begin{proof}
For any maximal ideal $\mathfrak{m}$ of $C$, it is easy to see that $(\text{tr}_{HS})_{\mathfrak{m}}$ is precisely the regular trace of $A_{\mathfrak{m}}$ as a finitely generated free $C_{\mathfrak{m}}$-module. It follows that $(A_{\mathfrak{m}},C_{\mathfrak{m}},(\text{tr}_{HS})_{\mathfrak{m}})$ is a Cayley-Hamilton algebra of degree $n=\text{rank}_{C_{\mathfrak{m}}}A_{\mathfrak{m}}$. For any $a\in A$, consider the $n$-th characteristic polynomial $p_{n,a}(t)\in C[t]$ defined by \eqref{eq-charpoly-CHalg} of $a$ with respect to $\text{tr}_{HS}$. 

 Note that $ p_{n,a}(a)/1=0\in C_{\mathfrak{m}}$, as $(A_{\mathfrak{m}},C_{\mathfrak{m}},(\text{tr}_{HS})_{\mathfrak{m}})$ is a Cayley-Hamilton algebra of degree $n$. Then the localization of the principal ideal generated by $p_{n,a}(a)$ is zero at any maximal ideal of $C$. This forces $p_{n,a}(a)=0$.   It follows similarly that $\text{tr}_{HS}(1)=n$.
\end{proof}


\subsection{Fiber algebras and indecomposable exact module categories}\label{sec-H/mHmocat-exact}
Let $\mathcal{C}$ be a finite tensor category over an algebraically closed field $\mathbbm{k}$. That is, $\mathcal{C}$ is a $\mathbbm{k}$-linear, abelian, rigid monoidal category with finite-dimensional Hom-spaces. Moreover, all objects have finite length, the unit object is irreducible, and there are only finitely many isomorphism classes of irreducible objects in $\mathcal{C}$. Throughout this subsection, $\mathbbm{k}$ is  an algebraically closed field of \textit{arbitrary} characteristic, the module category over $\mathcal{C}$ considered is by default a \textit{left $\mathcal{C}$-module category}. Suppose $\mathcal{M}$ is a finite module category over $\mathcal{C}$. Then the Grothendieck group $\text{Gr}(\mathcal{M})$ can be naturally regarded as a $\mathbb{Z}_{+}$-module over the Grothendieck ring $\text{Gr}(\mathcal{C})$ via the tensor action, with $\text{Irr}(\mathcal{M})$ as its fixed $\mathbb{Z}_{+}$-basis. 

Let $\mathcal{M}$ be a finite module category over $\mathcal{C}$. The module category $\mathcal{M}$ is called \textit{exact} if for any projective object $P$ in $\mathcal{C}$ and any object $M$ in $\mathcal{M}$, the object $P\otimes M$ is projective in $\mathcal{M}$. It is called \textit{indecomposable} if it cannot be equivalent to a direct sum of two nonzero $\mathcal{C}$-module categories. For additional background on module categories, readers are referred to \cite{MR3242743}.


We fix a pair $(H,C)$ of Hopf algebras over $\mathbbm{k}$ satisfying the hypothesis  (FinHopf). By \eqref{eq-righlefcomalh9}, any fiber algebra $H/\mathfrak{m}H$ is both a left and a right comodule algebra over the identity fiber algebra $H/\mathfrak{m}_{\overline{\varepsilon}}H$. It follows that $H/\mathfrak{m}H\text{-mod}$ can be regarded as both a left and a right module category over the tensor category $H/\mathfrak{m}_{\overline{\varepsilon}}H\text{-mod}$. More generally, according to \eqref{eq-coproun-algfibers}, for any $\chi,\psi\in G(C^{\circ})$, there is a bifunctor:
\begin{equation}\label{eq-bifunct-conv-fiberalgsten}
\otimes : H/\mathfrak{m}_{\chi}H\text{-mod} \times H/\mathfrak{m}_{\psi}H\text{-mod}\to H/\mathfrak{m}_{\chi\ast\psi}H\text{-mod}.
\end{equation}

The main goal of this subsection is to prove that for any maximal ideal $\mathfrak{m}$ of $C$, the category $H/\mathfrak{m}H\text{-mod}$, as a module category over the tensor category $H/\mathfrak{m}_{\overline{\varepsilon}}H\text{-mod}$, is exact and  indecomposable (Proposition \ref{thm-indexacmocat-fiberalg}). Before proving this, we need some further discussion.

Since $H$ is a Noetherian affine PI Hopf algebra, by \cite[Cor. 2]{MR2271358}, the antipode $S$ of $H$ is bijective. Consequently, the category of finite-dimensional $H$-modules is a tensor category with respect to the tensor product. For any finite-dimensional $H$-module $V$, the left dual $V^{*}$ and right dual $^{*}{V}$ of $V$ are the usual dual space of $V$, with actions of $H$ given by
$$(h\xi)(v)=\xi(S(h)v) \text{ and }(h\xi)(v)=\xi(S^{-1}(h)v)$$ respectively,
for all $h\in H, v\in V$ and $\xi$ in the dual space of $V$. If $V$ is an $H/\mathfrak{m}H$-module, then $V^{*}$ and $^{*}V$ are both $H/\mathfrak{m}^{-1}H$-modules. For any finite-dimensional $H$-modules $X,Y,Z$, one has the natural adjunction isomorphism (see also \cite[Prop. 2.10.8]{MR3242743})
\begin{equation}\label{eq-naadj-isous7}
\Phi:  \text{Hom}_{H}(X,Z\otimes Y^{*})\to \text{Hom}_{H}(X\otimes Y,Z),
\end{equation}
given by $\Phi(\theta)(x\otimes y)=(\text{id}_{Z}\otimes y)\theta(x)$. Now we are ready to prove the following.



\begin{proposition}
Let $(H,C)$ be a pair of Hopf algebras satisfying the hypothesis (FinHopf). Then for any $\mathfrak{m}\in \text{maxSpec}C$, $H/\mathfrak{m}H\text{-mod}$ is a finite indecomposable exact module category over the finite tensor category $H/\mathfrak{m}_{\overline{\varepsilon}}H\text{-mod}$.\label{thm-indexacmocat-fiberalg}
\end{proposition}
\begin{proof}
We first prove that $H/\mathfrak{m}H\text{-mod}$ is exact as a module category over $H/\mathfrak{m}_{\overline{\varepsilon}}H\text{-mod}$. By definition, we need to verify that for any projective $H/\mathfrak{m}_{\overline{\varepsilon}}H$-module $P$ and any $H/\mathfrak{m}H$-module $Y$, the tensor product $P\otimes Y$, considered as an $H/\mathfrak{m}H$-module via \eqref{eq-bifunct-conv-fiberalgsten}, is projective.

In this case, the adjunction formula \eqref{eq-naadj-isous7} induces a natural isomorphism
$$\text{Hom}_{H/\mathfrak{m}H}(P\otimes Y,-)\cong \text{Hom}_{H/\mathfrak{m}_{\overline{\varepsilon}}H}(P,-\otimes Y^{*}).$$
It follows that $P\otimes Y$ is a projective $H/\mathfrak{m}H$-module.

We then proceed to prove that $H/\mathfrak{m}H\text{-mod}$ is indecomposable as a module category over $H/\mathfrak{m}_{\overline{\varepsilon}}H\text{-mod}$. Suppose on the contrary that there exist nonzero module categories $\mathcal{M}_{1}$ and $\mathcal{M}_{2}$ over the tensor category $H/\mathfrak{m}_{\overline{\varepsilon}}H\text{-mod}$ such that,  $H/\mathfrak{m}H\text{-mod}$ is equivalent to $\mathcal{M}_{1}\oplus \mathcal{M}_{2}$ as module categories over $H/\mathfrak{m}_{\overline{\varepsilon}}H\text{-mod}$. Let $F:H/\mathfrak{m}H\text{-mod}\to \mathcal{M}_{1}\oplus \mathcal{M}_{2}$ be an equivalence of module categories over $H/\mathfrak{m}_{\overline{\varepsilon}}H\text{-mod}$. 
For any irreducible $H/\mathfrak{m}H$-module $W$, 
$F(W)$ is an irreducible object in $\mathcal{M}_{1}\oplus \mathcal{M}_{2}$. Thus $F(W)$ must be of the form $(X,0)$ or $(0,Y)$, where $X$ is irreducible in $\mathcal{M}_{1}$, $Y$ is irreducible in $\mathcal{M}_{2}$. Without loss of generality, suppose there is an irreducible $H/\mathfrak{m}H$-module $W_{0}$ such that $F(W_{0})=(X,0)$, where $X$ is an irreducible object in $\mathcal{M}_{1}$. Let $W$ be an irreducible $H/\mathfrak{m}H$-module. We claim that there exists an irreducible object $X_{W}$ in $\mathcal{M}_{1}$ such that $F(W)\cong (X_{W},0)$. It follows from \cite[Thm. 3.1 (a)]{mi2025lowest} that there is an irreducible $H/\mathfrak{m}_{\overline{\varepsilon}}H$-module $V$ such that $V\otimes W_{0} $ has a composition factor isomorphic to $W$. 
Hence $F(W)$ is a subquotient object of $F(V\otimes W_{0})\cong V\otimes F(W_{0})\cong (V\otimes X,0)$. Therefore, $F(W)\cong (X_{W},0)$ for some irreducible object $X_{W}$ in $\mathcal{M}_{1}$. The claim is proved.

Now, for any irreducible  object $Y$ in $\mathcal{M}_{2}$, there exists an irreducible $H/\mathfrak{m}H$-module $W^{\prime}$ such that $F(W^{\prime})\cong (0,Y)$ as $F$ is an equivalence. By the previous claim there exists an irreducible object $X_{W^{\prime}}$ in $\mathcal{M}_{1}$ such that $F(W^{\prime})\cong (X_{W^{\prime}},0)$. Hence $(X_{W^{\prime}},0)\cong (0,Y)$, which is a contradiction.
\end{proof}

Proposition \ref{thm-indexacmocat-fiberalg} indicates that classifying the indecomposable exact module categories over $H/\mathfrak{m}_{\overline{\varepsilon}}H\text{-mod}$ benefits the study of the representation theory of $H$.

It immediately follows from \cite[Prop. 7.7.2]{MR3242743} and Proposition \ref{thm-indexacmocat-fiberalg} that
\begin{corollary}\label{cor-GrirreZ+-hopf}
For any $\mathfrak{m}\in \text{maxSpec}C$, the Grothendieck group $\text{Gr}(H/\mathfrak{m}H)$ is an irreducible $\mathbb{Z}_{+}$-module over the Grothendieck ring $\text{Gr}(H/\mathfrak{m}_{\overline{\varepsilon}}H)$. 
\end{corollary}
\begin{remark}\label{rmk-irreZ+-appnewpf}
Corollary \ref{cor-GrirreZ+-hopf} enables us to capture much of the representation-theoretic information of fiber algebras. Suppose $H/\mathfrak{m}_{\overline{\varepsilon}}H$ is a basic algebra (for instance, the pair of Hopf algebras $(\mathcal{O}_{\epsilon}(G),C_{\epsilon}(G))$ given in Example \ref{eg-quantizecoornri}).
Then, the Grothendieck ring $\text{Gr}(H/\mathfrak{m}_{\overline{\varepsilon}}H)$ of the identity fiber algebra $H/\mathfrak{m}_{\overline{\varepsilon}}H$ is the group ring of the group of characters 
$$G_{0}\coloneqq G((H/\mathfrak{m}_{\overline{\varepsilon}}H)^{\circ}).$$
By the proof of \cite[Prop. 1.1 (i)]{MR1333750}, any irreducible $\mathbb{Z}_{+}$-module over $\mathbb{Z}G_{0}$ is of the form $\mathbb{Z}[G_{0}/N]$, where $N$ is some subgroup of $G_{0}$. In particular, for any $\mathfrak{m}\in \text{maxSpec}C$ the action of $G_{0}$ on $\text{Irr}(H/\mathfrak{m}H)$ given by tensor product is \textit{transitive}, and $|\text{Irr}(H/\mathfrak{m}H)|$ divides $|G_{0}|=|\text{Irr}(H/\mathfrak{m}_{\overline{\varepsilon}}H)|$. The preceding discussion provides a tensor-categorical proof of \cite[Prop. III.4.11 and Cor. III.4.11]{MR1898492}. Note also that the transitivity of the group action of $G_{0}$ on $\text{Irr}(H/\mathfrak{m}H)$ implies that any two irreducible $H/\mathfrak{m}H$-modules have the same $\mathbbm{k}$-dimension. This observation, however, fails if we only assume that the identity fiber algebra has the Chevalley property (see \S \ref{sec-chevloc} for the definition). In fact, the Hopf algebra $H$ given in Example \ref{eg-Chevnotsembas-dim16} is a finite-dimensional Hopf algebra with the Chevalley property that is not basic. If we take $C=\mathbbm{k}$, then the identity fiber algebra $H/\mathfrak{m}_{\overline{\varepsilon}}H\cong H$ has two irreducible modules with different $\mathbbm{k}$-dimensions.
\end{remark}
We conclude this subsection with the following corollary.
\begin{corollary}
If $H/\mathfrak{m}_{\overline{\varepsilon}}H$ is semisimple, then all the fiber algebras $H/\mathfrak{m}H$ are semisimple.\label{cor-semisimall1}
\end{corollary}
\begin{proof}
Since the trivial module $_{H/\mathfrak{m}_{\overline{\varepsilon}}H}\mathbbm{k}$ is projective, it follows from Proposition \ref{thm-indexacmocat-fiberalg} that any finite-dimensional $H/\mathfrak{m}H$-module is projective. In particular, all irreducible $H/\mathfrak{m}H$-modules are projective.  It follows that any $H/\mathfrak{m}H$-module of finite length is completely reducible. Therefore, 
$H/\mathfrak{m}H$ is completely reducible as a module over itself, which implies that $H/\mathfrak{m}H$ is semisimple.
\end{proof}
\begin{remark}
There exist pairs $(H,C)$ of Hopf algebras satisfying the hypothesis (FinHopf) such that the identity fiber algebra $H/\mathfrak{m}_{\overline{\varepsilon}}H$ is not semisimple, while all fiber algebras $H/\mathfrak{m}H$ at $\mathfrak{m}\neq \mathfrak{m}_{\overline{\varepsilon}}\in \text{maxSpec}C$ are semisimple, as given in Example \ref{eg-infiTaft}. 
\end{remark}
\begin{remark}
For any pair $(H,C)$ of Hopf algebras over an algebraically closed field $\mathbbm{k}$ of characteristic zero satisfying (FinHopf), under the assumption that the identity fiber algebra $H/\mathfrak{m}_{\overline{\varepsilon}}H$ is semisimple, Brown-Couto \cite[Thm. 2.9 (8)]{MR4201485} proved that $H$ is AS-regular.
\end{remark}

\subsection{Orbits of the winding automorphism group action}\label{subsec-orbiwinda}

Let $(H,C)$ be a pair of Hopf algebras satisfying the hypothesis (FinHopf) over an algebraically closed field  $\mathbbm{k}$ of \textit{arbitrary} characteristic. The \textit{left winding automorphism} $W_{l}(\chi)$ and the \textit{right winding automorphism} $W_{r}(\chi)$ of $H$ associated to a character $\chi \in G(H^{\circ})$  are defined respectively as:
$$W_{l}(\chi)(h)\coloneqq \sum_{(h)}\chi(h_{(1)})h_{(2)} \text{ and\ } W_{r}(\chi)(h)\coloneqq \sum_{(h)}h_{(1)}\chi(h_{(2)}), \forall h\in H.$$

Winding automorphisms are crucial in the study of Noetherian Hopf algebras. Since Noetherian affine PI Hopf algebras are AS-Gorenstein \cite[Thm. 0.2 (1)]{MR1938745}, it follows that $H$ is an AS-Gorenstein Hopf algebra. By \cite[Prop. 4.5 (a)]{MR2437632}, there is a character $\chi_{0}\in G(H^{\circ})$, defined by the \textit{homological integral} of $H$ in the sense of Lu-Wu-Zhang \cite{MR2320655}, such that $W_{l}(\chi_{0})S^{2}$ is a \textit{Nakayama automorphism} of $H$, which is unique up to an inner automorphism.


Clearly, for any $\chi\in G(H^{\circ})$, one has $W_{l}(\chi)(C)=C$ and $W_{r}(\chi)(C)=C$. It follows that both $W_{l}(\chi)$ and $W_{r}(\chi)$ act on $\text{maxSpec}C$ as homeomorphisms. Moreover, $W_{l}(G(H^{\circ}))$ and $W_{r}(G(H^{\circ}))$ are both subgroups of the algebra automorphism group of $H$, referred to as the \textit{left winding automorphism group} and the \textit{right winding automorphism group} of $H$, respectively.

We consider the cosets of $\text{maxSpec}C$ with respect to the subgroup
\begin{equation}\label{subgroup-I}
 \mathcal{I} \coloneqq \{\mathfrak{m}\in \text{maxSpec}C \mid  H/\mathfrak{m}H\text{ has a }1\text{-dimensional module}\}.   
\end{equation}

The main result of this section is to prove that for any $\mathfrak{m}\in \text{maxSpec}C$, the \textit{left} coset of $\mathcal{I}$ in $\text{maxSpec}C$ which contains $\mathfrak{m}$ is exactly the orbit of $\mathfrak{m}$ of the \textit{right} winding automorphism group action (Theorem \ref{thm-leftcossta-isofiber-or}). 

The importance of the subgroup $\mathcal{I}$ is based on the fact that the points in $\mathcal{I}$ are exactly those maximal ideals of $C$ for which the corresponding fiber algebras are isomorphic to the identity fiber algebra $H/\mathfrak{m}_{\overline{\varepsilon}}H$. Equivalently, $\mathcal{I}$ is the orbit of $\mathfrak{m}_{\overline{\varepsilon}}$ of the left (or right) winding automorphism group action. In fact, for any maximal ideal $\mathfrak{m}$ of $C$, it is proved in \cite[Thm. 4.3 (a)]{mi2025lowest} that the following are equivalent, where an additional requirement that the identity fiber algebra $H/\mathfrak{m}_{\overline{\varepsilon}}H$ is a basic algebra is added.  
\begin{itemize}
\item[(1)] $\mathfrak{m}$ belongs to the subgroup $\mathcal{I}$;
\item[(2)] $H/\mathfrak{m}H$ and $H/\mathfrak{m}_{\overline{\varepsilon}}H$ are isomorphic as $\mathbbm{k}$-algebras;
\item[(3)] $\mathfrak{m}$ belongs to the orbit of $\mathfrak{m}_{\overline{\varepsilon}}$ of the right winding automorphism group action;
\item[(4)] $\mathfrak{m}$ belongs to the orbit of $\mathfrak{m}_{\overline{\varepsilon}}$ of the left winding automorphism group action.
\end{itemize}

Now, we proceed to prove the main result of this subsection.
\begin{theorem}\label{thm-leftcossta-isofiber-or}
Let $(H,C)$ be a pair of Hopf algebras satisfying the hypothesis (FinHopf), and let $\mathcal{I}$ be the subgroup of $\text{maxSpec}C$ defined in \eqref{subgroup-I}. Then for any  maximal ideal $\mathfrak{m}$ of $C$, the left coset of $\mathcal{I}$ in $\text{maxSpec}C$ containing $\mathfrak{m}$ coincides with the orbit of $\mathfrak{m}$ of the right winding automorphism group action.
\end{theorem}

\begin{remark}
Theorem \ref{thm-leftcossta-isofiber-or} shows that for maximal ideas $\mathfrak{m}$ and $\mathfrak{n}$ of $C$ belonging to the same left coset of $\mathcal{I}$ in $\text{maxSpec}C$, $H/\mathfrak{m}H$ and $H/\mathfrak{n}H$ are isomorphic as $\mathbbm{k}$-algebras. In general, there exist pairs $(H,C)$ satisfying (FinHopf) such that the fiber algebras at maximal ideals from any two different left cosets of $\mathcal{I}$ are non-isomorphic. See Example \ref{eg-Chevnotsembas-dim16}. We remark that Theorem \ref{thm-leftcossta-isofiber-or} can also be regarded as an improvement of \cite[Thm. 4.3 (a)]{mi2025lowest}.
\end{remark}
\begin{proof}[Proof of Theorem \upshape\ref{thm-leftcossta-isofiber-or}]
Let $\mathfrak{m}_{\varphi}$ and $\mathfrak{m}_{\psi}$ be maximal ideals of $C$ in the same left coset of $\mathcal{I}$ in $\text{maxSpec}C$. Then there exists $\mathfrak{m}_{\chi}\in \mathcal{I}$ such that $\mathfrak{m}_{\psi}= \mathfrak{m}_{\varphi}\ast \mathfrak{m}_{\chi}$. It follows that $\psi= \varphi\ast \chi$. According to the definition of $\mathcal{I}$, the fiber algebra $H/\mathfrak{m}_{\chi}H$ has a $1$-dimensional representation $\rho:H/\mathfrak{m}_{\chi}H\to\mathbbm{k}$. For $h\in H$, let $\overline{h}=h+\mathfrak{m}_{\chi}H\in H/\mathfrak{m}_{\chi}H$. Define $\widehat{\chi}:H\to\mathbbm{k}$ by $h\mapsto \rho(\overline{h})$. Then $\widehat{\chi}\in G(H^{\circ})$ and $\widehat{\chi}|_{C}=\chi$.

 For any $c\in \mathfrak{m}_{\psi}=\text{Ker}\psi$, one has 
$$0=\psi(c)=(\varphi\ast\chi)(c)=\sum_{(c)}\varphi(c_{(1)})\chi(c_{(2)})=\varphi\left(\sum_{(c)}c_{(1)}\chi(c_{(2)})\right)=\varphi(W_{r}(\hat{\chi})(c)).$$
Thus $W_{r}(\widehat{\chi})(\mathfrak{m}_{\psi})=\mathfrak{m}_{\varphi}$. Therefore, for any maximal ideal $\mathfrak{m}$ of $C$, the left coset $\mathfrak{m}\ast \mathcal{I}$ is a subset of the orbit $W_{r}(G(H^{\circ}))(\mathfrak{m})$. Now, we show that $W_{r}(G(H^{\circ}))(\mathfrak{m})\subseteq \mathfrak{m}\ast \mathcal{I}$ to complete the proof. Suppose $\mathfrak{m}_{\psi}\in W_{r}(G(H^{\circ}))(\mathfrak{m})$ and let $\mathfrak{m}=\mathfrak{m}_{\varphi}$. Then there exist $\widehat{\chi}\in G(H^{\circ})$ such that $\mathfrak{m}_{\psi}=W_{r}(\widehat{\chi})(\mathfrak{m}_{\varphi})$. Set $\chi\coloneqq \widehat{\chi}|_{C}\in G(C^{\circ})$. We claim that $\mathfrak{m}_{\varphi}=\mathfrak{m}_{\psi}\ast \mathfrak{m}_{\chi}$.

Take $c\in \mathfrak{m}_{\varphi}$. Then we have $\psi(W_{r}(\widehat{\chi})(c))=0$. Thus
$$0=\psi\left(\sum_{(c)}c_{(1)}\widehat{\chi}(c_{(2)})\right)=\psi\left(\sum_{(c)}c_{(1)}\chi(c_{(2)})\right)=(\psi\ast \chi)(c).$$
It follows that $\mathfrak{m}_{\varphi}\subseteq \mathfrak{m}_{\psi\ast \chi}$, which implies that $\mathfrak{m}_{\varphi}=\mathfrak{m}_{\psi}\ast \mathfrak{m}_{\chi}$.
\end{proof}
\begin{remark}\label{rmk-ricoscas-wnorsub1}
A similar argument shows that the \textit{right} cosets of $\mathcal{I}$ in $\text{maxSpec}C$ correspond to the orbits of the \textit{left} winding automorphism group action. Moreover, it is not hard to prove, by Theorem \ref{thm-leftcossta-isofiber-or}, that when $\mathcal{I}$ is a normal subgroup of $\text{maxSpec}C$, for any maximal ideal $\mathfrak{m}$ of $C$, the orbit of $\mathfrak{m}$ of the left winding automorphism group action coincides with the orbit of the right winding automorphism group action. We leave the details to the reader.
\end{remark}

\section{\texorpdfstring{$\overline{\varepsilon}$-Chevalley locus}{epsilon Chevalley locus} and discriminant ideals}\label{sec-Chevloc-disccon}
In this section, $\mathbbm{k}$ is assumed to be an algebraically closed field of characteristic zero. 

\subsection{The Chevalley property and the \texorpdfstring{$\overline{\varepsilon}$-Chevalley locus}{epsilon Chevalley locus}}\label{sec-chevloc}
Let $(H,C)$ be a pair of Hopf algebras satisfying the hypothesis (FinHopf). In this subsection, after reviewing the background of the Chevalley property of Hopf algebras, we introduce the concept of the \textit{$\overline{\varepsilon}$-Chevalley locus} (Definition \ref{def-Chev-locus}) for the pair $(H,C)$, which captures those fiber algebras whose completely reducible representations remain completely reducible under the tensor action with the irreducible representations of the identity fiber algebra. We will prove that the $\overline{\varepsilon}$-Chevalley locus of $(H,C)$ is the entire maximal spectrum $\text{maxSpec}C$ of $C$ if and only if the identity fiber algebra $H/\mathfrak{m}_{\overline{\varepsilon}}H$ has the Chevalley property in Proposition \ref{prop-epiChev-mepiC}.


We first recall the definition of the Chevalley property of Hopf algebras. 


A Hopf algebra $H$ is said to have the \textit{Chevalley property} \cite{MR1852304} if the tensor product of any two finite-dimensional irreducible $H$-modules is completely reducible (see also \cite[Def. 4.12.3]{MR3242743} for the definition of the Chevalley property for tensor categories). This terminology is named after a classical theorem of Chevalley, which states that for any group (not necessary finite), the tensor product of any two finite-dimensional irreducible representations over $\mathbb{C}$ is completely reducible.

Suppose $H$ is a finite-dimensional Hopf algebra over $\mathbbm{k}$. 
It was proved in \cite[Prop. 4.2]{MR1852304} that the following conditions are equivalent: 
\begin{itemize}
\item[(1)] $H$ has the Chevalley property.
\item[(2)] The Jacobson radical of $H$ is a Hopf ideal of $H$.
\item[(3)] The coradical of $H^{\circ}$ is a Hopf subalgebra of $H^{\circ}$.
\end{itemize}

If a finite-dimensional Hopf algebra $H$ is basic (or equivalently, $H^{\circ}$ is pointed), or semisimple, then it has the Chevalley property. The converse is not true as the following example shows.

\begin{example}
Let $H$ be the algebra over $\mathbbm{k}$ generated by $b,c,x,y$, subject to the following relations \cite[\S 5]{MR2037722}:\label{eg-Chevnotsembas-dim16}
$$ b^{2}=1,\ c^{2}=1,  \ x^{2}=\frac{1}{2}(1+c+b-cb),$$
$$bc=cb,\ xc=bx,\ xb=cx,$$
$$y^{2}=0,\ yc=-cy,\ yb=-by,yx=\sqrt{-1}cxy.$$
The coalgebra structure of $H$ is defined by
$$\Delta(b)=b\otimes b,\ \varepsilon(b)=1,$$
$$\Delta(c)=c\otimes c,\ \varepsilon(c)=1,$$
$$\Delta(x)=\frac{1}{2}(x\otimes x +bx\otimes x+x\otimes cx-bx\otimes cx),\ \varepsilon(x)=1,$$
$$\Delta(y)=y\otimes 1+c\otimes y,\ \varepsilon(y)=0. $$
The antipode $S$ of $H$ is defined by
$$S(b)=b=b^{-1},\ S(c)=c=c^{-1},\ S(x)=x,\text{ and}\ S(y)=-c^{-1}y=-cy.$$
Then $H$ is a Hopf algebra of dimension $16$, and $H$ has four $1$-dimensional representations:
$$\rho_{1}(b)=\rho_{1}(c)=1,\ \rho_{1}(x)=1,\ \rho_{1}(y)=0;$$
$$\rho_{2}(b)=\rho_{2}(c)=1,\ \rho_{2}(x)=-1,\ \rho_{2}(y)=0;$$
$$\rho_{3}(b)=\rho_{3}(c)=-1,\ \rho_{3}(x)=\sqrt{-1},\ \rho_{3}(y)=0;$$
$$\rho_{4}(b)=\rho_{4}(c)=-1,\ \rho_{4}(x)=-\sqrt{-1},\ \rho_{4}(y)=0;$$
and one $2$-dimensional irreducible representation
$$\rho_{5}(b)=\begin{pmatrix}
-1&0\\
0& 1
\end{pmatrix},\ 
\rho_{5}(c)=\begin{pmatrix}
1&0\\
0& -1
\end{pmatrix},\
\rho_{5}(x)=\begin{pmatrix}
0&1\\
1& 0
\end{pmatrix},\
\rho_{5}(y)=\begin{pmatrix}
0&0\\
0& 0
\end{pmatrix}.
$$
It can be  verified directly that the Jacobson radical of $H$ is the ideal generated by $y$, which is a Hopf ideal of $H$.  It follows that $H$ is a non-basic, non-semisimple Hopf algebra with the Chevalley property. The semisimple Hopf algebra $H/(y)$ is known as the \textit{Kac-Paljutkin's Hopf algebra} \cite{MR208401}. Moreover, $H$ has a central Hopf subalgebra $C=\mathbbm{k}\oplus \mathbbm{k}bc\cong \mathbbm{k}\mathbb{Z}_{2}$. Thus, $C$ has two maximal ideals: $\mathfrak{m}_{\overline{\varepsilon}}=\{\alpha-\alpha bc\in C\mid \alpha\in \mathbbm{k}\}, \text{and\ } \mathfrak{m}=\{\alpha+\alpha bc\in C\mid\alpha\in \mathbbm{k}\}.$ In this case, both fiber algebras $H/\mathfrak{m}_{\overline{\varepsilon}}H$ and $H/\mathfrak{m}H$ are non-semisimple of dimension $8$. For the pair $(H,C)$, the subgroup $\mathcal{I}$ defined in \eqref{subgroup-I} is the trivial group $\{\mathfrak{m}_{\overline{\varepsilon}}\}$. Hence, the two left cosets of $\mathcal{I}$ correspond to distinct isomorphism classes of fiber algebras as singleton sets. It is also clear that $(H,C,\text{tr}_{\text{reg}})$ is a Cayley-Hamilton Hopf algebra of degree $8$.

\end{example}

Next we prove that for any  pair $(H,C)$ satisfying the hypothesis (FinHopf) the semisimplicity of the identity fiber algebra $H/\mathfrak{m}_{\overline{\varepsilon}}H$ implies that $H$ has the Chevalley property. Note that Example \ref{eg-Chevnotsembas-dim16} also provides a pair $(H,C)$ of Hopf algebras that satisfies the hypothesis (FinHopf) such that $H$ has the Chevalley property but the identity fiber algebra $H/\mathfrak{m}_{\overline{\varepsilon}}H$ is not semisimple. 

\begin{proposition}
Let $(H,C)$ be a pair of Hopf algebras satisfying the hypothesis (FinHopf). If $H/\mathfrak{m}_{\overline{\varepsilon}}H$ is semisimple, then $H$ has the Chevalley property.\label{prop-H/mepiHsem-chev}
\end{proposition}
\begin{proof}
By Corollary \ref{cor-semisimall1}, all the fiber algebras $H/\mathfrak{m}H$ are semisimple. For any $\mathfrak{m}, \mathfrak{n} \in \text{maxSpec}C$, if $V$ and $W$ are  $H/\mathfrak{m}H$-module  and $H/\mathfrak{n}H$-module respectively, by \eqref{eq-coproun-algfibers}, $V\otimes W$ can be viewed as an $H/(\mathfrak{m}\ast\mathfrak{n})H$-module. It follows that $V\otimes W$ is completely reducible as an $H/(\mathfrak{m}\ast\mathfrak{n})H$-module. Consequently, $V\otimes W$ is also a completely reducible $H$-module. The conclusion follows.
\end{proof}
An immediate consequence of Proposition \ref{prop-H/mepiHsem-chev} is the following, as $\text{char}\mathbbm{k}=0$ is assumed.
\begin{corollary}
Let $(H,C)$ be a pair of Hopf algebras satisfying the hypothesis (FinHopf). If $H$ is involutory, that is, $S^{2}=\text{id}$, then $H$ has the Chevalley property.\label{cor-invrHop-chev}
\end{corollary}
\begin{proof}
Note that $H/\mathfrak{m}_{\overline{\varepsilon}}H$ is a finite-dimensional involutory Hopf algebra over a field of characteristic zero. It follows from a classical result of Larson and Radford \cite{MR957441,MR926744} that $H/\mathfrak{m}_{\overline{\varepsilon}}H$ is semisimple. The conclusion follows from Proposition \ref{prop-H/mepiHsem-chev}.
\end{proof}
Let $(H,C)$ be a pair of Hopf algebras satisfying the hypothesis (FinHopf). If  $H$ is commutative or cocommutative, then $H$ is involutory, and by Corollary \ref{cor-invrHop-chev}, $H$ has the Chevalley property.
\begin{example}\label{eg-grpalg-cenex1}
Let $G$ be a finitely generated group with a central subgroup $N$ of finite index. Let $\{g_{1},...,g_{m}\}$ be a set of right coset representatives of $G$ with respect to $N$, where $m=[G:N]$. Set $H=\mathbbm{k}G$ and $C=\mathbbm{k}N$. It follows that $N$ is a finitely generated abelian group. Thus, $G$ is a central extension of the finite group $G/N$ by a finitely generated abelian group $N$. Then $H$ is an affine Hopf algebra and $H$ is finitely generated as a module over the central Hopf subalgebra $C$. It is also clear that the pair $(H,C)$ satisfies the hypothesis (FinHopf) and that $H$ is a free module of rank $m$, with basis $\{g_{1},...,g_{m}\}$. Consequently, the $\mathbbm{k}$-dimension of each fiber algebra in $\mathcal{B}_{H}=\{H/\mathfrak{m}H\mid \mathfrak{m}\in \text{maxSpec}C\}$ is $m$. Since $H$ is involutory, the identity fiber algebra $H/\mathfrak{m}_{\overline{\varepsilon}}H$ is clearly semisimple, and $H$ has the Chevalley property by Corollary \ref{cor-invrHop-chev}. In fact, the natural projection $H=\mathbbm{k}G\to \mathbbm{k}[G/N]$ induces an algebra isomorphism $H/\mathfrak{m}_{\overline{\varepsilon}}H\cong \mathbbm{k}[G/N]$. By Maschke's Theorem, the group algebra $\mathbbm{k}[G/N]$ is semisimple. It is clear that $\mathbbm{k}[G/N]$ is basic if and only if $G/N$ is abelian.
\end{example}

Next we define the $\overline{\varepsilon}$-Chevalley locus for a pair $(H,C)$ satisfying the hypothesis (FinHopf).
\begin{definition}
Let $(H,C)$ be a pair of Hopf algebras satisfying the hypothesis (FinHopf). 
\begin{align*}
\text{Chev}_{\overline{\varepsilon}}(H,C)\coloneqq \{\mathfrak{m}\in \text{maxSpec}C\mid &V\otimes W  \text{\ is a completely reducible\ } H/\mathfrak{m}H\text{-module for all\ }\\
&[V]\in \text{Irr}(H/\mathfrak{m}_{\overline{\varepsilon}}H)\text{ and\ }[W]\in \text{Irr}(H/\mathfrak{m}H)\}
\end{align*}
is deined to be the \textit{$\overline{\varepsilon}$-Chevalley locus} of the pair $(H,C)$.\label{def-Chev-locus}
\end{definition}

\begin{remark}In general, the $\overline{\varepsilon}$-Chevalley locus may be empty for a pair $(H,C)$ of Hopf algebras satisfying the hypothesis (FinHopf). Let $C=\mathbbm{k}$ for a finite-dimensional Hopf algebra $H$. Then $H$ has the Chevalley property if and only if $\text{Chev}_{\overline{\varepsilon}}(H,C)=\text{maxSpec}C=\{\mathfrak{m}_{\overline{\varepsilon}}\}=\{0\}$. Thanks to \cite[Cor. III.4.10]{MR1898492}, for a pair $(H,C)$ satisfying the hypothesis (FinHopf) such that $H$ is a prime ring, its $\overline{\varepsilon}$-Chevalley locus is a dense subset of the  affine variety $\text{maxSpec}C$. Indeed, $\mathfrak{m}\in  \text{Chev}_{\overline{\varepsilon}}(H,C)$ for any $\mathfrak{m}\in \text{maxSpec}C$ such that $H/\mathfrak{m}H$ is a semisimple algebra.
\end{remark}
\begin{remark}
By definition, the $\overline{\varepsilon}$-Chevalley locus of the pair $(H,C)$ clearly depends on $C$. There are Hopf algebras $H$ with different central Hopf subalgebras $C_{1}$ and $C_{2}$ such that both $(H,C_{1})$ and $(H,C_{2})$ satisfy the hypothesis (FinHopf), but $\text{Chev}_{\overline{\varepsilon}}(H,C_{1})=\text{maxSpec}C_{1}$ and $\text{Chev}_{\overline{\varepsilon}}(H,C_{2}) \subsetneqq \text{maxSpec}C_{2}$, as  Example \ref{eg-cheloc=maxSp-notChprop} shows. 
\end{remark}
 For any pair $(H,C)$ of Hopf algebras satisfying the hypothesis (FinHopf), it is clear that $\text{Chev}_{\overline{\varepsilon}}(H,C)=\text{maxSpec}C$ if $H$ has the Chevalley property. The following example shows that the converse is not true.
\begin{example} \label{eg-cheloc=maxSp-notChprop}
Consider the Hopf algebra $H=H(\mu)$ where $\mu\in \mathbbm{k}^{\times}$ \cite{MR4396654}.  As a $\mathbbm{k}$-algebra, $H$ is generated by $x,g$ subject to the relations 
$$g^{4}=1,\ xg=-gx,\ x^{2}=\mu(1-g^{2})/2.$$
The coalgebra structure of $H$ is defined by $\Delta(g)=g\otimes g,\Delta(x)=x\otimes 1+g\otimes x,\varepsilon(g)=1$ and $\varepsilon(x)=0$. The antipode $S$ of $H$ is given by $S(g)=g^{-1}$ and $S(x)=-g^{-1}x$. Then $H$ is a Hopf algebra of dimension $8$, with a central Hopf subalgebra $C=\mathbbm{k}[g^{2}]$. It is direct to verify that $H$ is a free $C$-module with basis $\{1,x,g,gx\}$ and $C$ possesses two maximal ideals $\mathfrak{m}_{\overline{\varepsilon}}=(g^{2}-1)$ and $\mathfrak{m}=(g^{2}+1)$. Clearly, the identity fiber algebra $H/\mathfrak{m}_{\overline{\varepsilon}}H$ is the  $4$-dimensional Taft algebra $H_{2}(-1)=\mathbbm{k}\langle g,x\rangle/(g^{2}-1,xg+gx,x^{2})$, which is basic. One also has an algebra homomorphism $\zeta:H/\mathfrak{m}H=H/(g^{2}+1)H\to\text{M}_{2}(\mathbbm{k})$ defined by
$$\zeta(g)=
\begin{pmatrix}
\sqrt{-1}&0\\
0&-\sqrt{-1}
\end{pmatrix},
\zeta(x)=
\begin{pmatrix}
0&\mu\\
1&0
\end{pmatrix}.
$$
Since $\mu\in \mathbbm{k}^{\times}$, it can be easily checked that $\zeta$ is surjective. It follows from $\dim_{\mathbbm{k}}H/\mathfrak{m}H=4=\dim_{\mathbbm{k}}\text{M}_{2}(\mathbbm{k})$ that $\zeta$ is bijective. Hence, $H/\mathfrak{m}H\cong \text{M}_{2}(\mathbbm{k})$ as $\mathbbm{k}$-algebras. Let $V$ be an irreducible $H/\mathfrak{m}H$-module, which is unique up to isomorphism. Since $\mathfrak{m}=\mathfrak{m}^{-1}$, the left dual $V^{*}$ is an irreducible $H/\mathfrak{m}H$-module. Hence $V\cong V^{*}$ as $H/\mathfrak{m}H$-modules. We claim that $V\otimes V$ is not a completely reducible $H$-module. Suppose on the contrary that $V\otimes V$ is a completely reducible $H$-module. It follows from $V\cong V^{*}$ that $V\otimes V^{*}$ is a completely reducible $H/\mathfrak{m}_{\overline{\varepsilon}}H$-module. Hence 
the trivial module $_{H/\mathfrak{m}_{\overline{\varepsilon}}H}\mathbbm{k}$ is a direct summand of $V\otimes V^{*}$.
By the adjunction formula \eqref{eq-naadj-isous7}, one has the natural isomorphisms \begin{equation}\label{eq-eg8ch-isof}
\text{Hom}_{H/\mathfrak{m}_{\overline{\varepsilon}}H}(V\otimes V^{*},-)\cong \text{Hom}_{H/\mathfrak{m}H}(V,-\otimes V^{**})\cong \text{Hom}_{H/\mathfrak{m}H}(V,-\otimes V).
\end{equation}
Since $V$ is a projective $H/\mathfrak{m}H$-module, it follows from \eqref{eq-eg8ch-isof}
that $V\otimes V^{*}$ is a projective $H/\mathfrak{m}_{\overline{\varepsilon}}H$-module. Hence the trivial module $_{H/\mathfrak{m}_{\overline{\varepsilon}}H}\mathbbm{k}$ is projective, which implies that the identity fiber algebra $H/\mathfrak{m}_{\overline{\varepsilon}}H$ is semisimple. This contradicts the fact that $H/\mathfrak{m}_{\overline{\varepsilon}}H\cong H_{2}(-1)$ is not semisimple. Thus $H$ does not have the Chevalley property. However, one has $\text{Chev}_{\overline{\varepsilon}}(H,C)=\text{maxSpec}C$ by \cite[Thm. 3.1 (b)]{mi2025lowest}. By Example \ref{eg-regulartr}, $(H,C,\text{tr}_{\text{reg}})$ is a Cayley-Hamilton Hopf algebra of degree $4$.
\end{example}

So, $\text{Chev}_{\overline{\varepsilon}}(H,C)=\text{maxSpec}C$
is not sufficient to ensure that $H$ has the Chevalley property. It is, however, equivalent to $H/\mathfrak{m}_{\overline{\varepsilon}}H$ having the Chevalley property.

\begin{proposition}
Let $(H,C)$ be a pair of Hopf algebras satisfying the hypothesis (FinHopf). Then $\text{Chev}_{\overline{\varepsilon}}(H,C)=\text{maxSpec}C$ if and only if $H/\mathfrak{m}_{\overline{\varepsilon}}H$ has the Chevalley property.\label{prop-epiChev-mepiC}
\end{proposition}
\begin{proof}
We need only to prove the ``if part". Since the Jacobson radical $\text{Jac}(H/\mathfrak{m}_{\overline{\varepsilon}}H)$ is a Hopf ideal of the identity fiber algebra $H/\mathfrak{m}_{\overline{\varepsilon}}H$, $R\coloneqq (H/\mathfrak{m}_{\overline{\varepsilon}}H)/\text{Jac}(H/\mathfrak{m}_{\overline{\varepsilon}}H)$ is a finite-dimensional semisimple Hopf algebra. For each $\mathfrak{m}\in \text{maxSpec}C$, the left $H/\mathfrak{m}_{\overline{\varepsilon}}H$-comodule algebra structure on $H/\mathfrak{m}H$ given by \eqref{eq-righlefcomalh9} induces a left $R$-comodule algebra structure on $H/\mathfrak{m}H$. Equivalently, $H/\mathfrak{m}H$ is a right $R^{\circ}$-module algebra. By \cite[Thm. 2.1]{MR1995055} or \cite[Thm. 0.3]{MR3359721}, $\text{Jac}(H/\mathfrak{m}H)$ is stable under the Hopf action of $R^{\circ}$. It follows that $\text{Jac}(H/\mathfrak{m}H)$ is a left $R$-subcomodule of $H/\mathfrak{m}H$. Hence, $\overline{\Delta}(\text{Jac}(H/\mathfrak{m}H))\subseteq  H/\mathfrak{m}_{\overline{\varepsilon}}H\otimes \text{Jac}(H/\mathfrak{m}H)+\text{Jac}(H/\mathfrak{m}_{\overline{\varepsilon}}H)\otimes H/\mathfrak{m}H$. This implies that $\mathfrak{m}\in \text{Chev}_{\overline{\varepsilon}}(H,C)$.
\end{proof}
\begin{remark}\label{rmk-H/mepHbasic-Chev}
In particular, $\text{Chev}_{\overline{\varepsilon}}(H,C)=\text{maxSpec}C$ if the identity fiber algebra $H/\mathfrak{m}_{\overline{\varepsilon}}H$ is either semisimple or basic. Let $(\mathcal{O}_{\epsilon}(G),C_{\epsilon}(G))$ be the pair of Hopf algebras as described in Example \ref{eg-quantizecoornri}. Then its identity fiber algebra $\mathcal{O}_{\epsilon}(G)/\mathfrak{m}_{\overline{\varepsilon}}\mathcal{O}_{\epsilon}(G)$ is basic. 
\end{remark}
\begin{remark}
Combining Example \ref{eg-cheloc=maxSp-notChprop} and Proposition \ref{prop-epiChev-mepiC}, whether the identity fiber algebra of a pair $(H,C)$ of Hopf algebras satisfying the hypothesis (FinHopf) has the Chevalley property depends on the choice of the central Hopf subalgebra $C$.
\end{remark}

To sum up, for any pair $(H,C)$ of Hopf algebras satisfying the hypothesis (FinHopf), 
\begin{align*}
H\text{ has the Chevalley property}\Rightarrow &\text{Chev}_{\overline{\varepsilon}}(H,C)=\text{maxSpec}C\\
\Leftrightarrow & H/\mathfrak{m}_{\overline{\varepsilon}}H\text{ has the Chevalley property.}
\end{align*}
In general, it is not easy to check whether a Hopf algebra (even if the Hopf algebra in question is finite-dimensional) has the Chevalley property by definition. 
In \S \ref{sec-Chev-alllrti2}, discriminant ideals (see Definition \ref{def-discriid}) are used to test whether an affine Hopf algebra admitting a large central Hopf subalgebra has the Chevalley property, see Theorem \ref{thm-CHHopf-Chev-dis}.



\subsection{The behavior of the square dimension function}\label{sec-behSqfun-Chev}
Let $A$ be an affine module-finite $C$-algebra. Consider the \textit{square dimension function} \cite{mi2025lowest}:
\begin{equation}\label{eq-Sd-function}
\text{Sd}:\text{maxSpec}C\to \mathbb{N}, \mathfrak{m}\mapsto \sum_{[V]\in \text{Irr}(A/\mathfrak{m}A)}(\dim_{\mathbbm{k}}V)^{2}.
\end{equation}
Suppose that $(H,C)$  is a pair of Hopf algebras satisfying the hypothesis (FinHopf) such that the identity fiber algebra $H/\mathfrak{m}_{\overline{\varepsilon}}H$ has the Chevalley property.
We investigate the behavior of the square dimension function of the pair $(H,C)$ on the maximal spectrum $\text{maxSpec}C$.

The following is the main result in this subsection, which is used to prove Theorem \ref{intro-thmB} (Theorem \ref{thm-Chein-VDneth-cwao}) in the introduction.
\begin{theorem}
Let $(H,C)$ be a pair of Hopf algebras satisfying the hypothesis (FinHopf) such that the identity fiber algebra $H/\mathfrak{m}_{\overline{\varepsilon}}H$ has the Chevalley property.\label{thm-sdchevpo-geqmepi} Set $\text{Irr}(H/\mathfrak{m}_{\overline{\varepsilon}}H)=\{[V_{1}],...,[V_{m}]\}$. Then for any $\mathfrak{m}\in \text{maxSpec}C$, 
\begin{itemize}
 \item[(1)] $\text{Sd}(\mathfrak{m})\geq \text{Sd}(\mathfrak{m}_{\overline{\varepsilon}})\geq |\text{Irr}(H/\mathfrak{m}H)|$.
 \item[(2)] The following are equivalent:
 \begin{itemize}
 \item[(i)]$\text{Sd}(\mathfrak{m})=\text{Sd}(\mathfrak{m}_{\overline{\varepsilon}})$;
 \item[(ii)] For any irreducible $H/\mathfrak{m}H$-module $W$, as $H$-modules, 
 $$W\otimes W^{*}\cong \bigoplus_{i=1}^{m}V_{i}^{\oplus \dim_{\mathbbm{k}}\text{Hom}_{H}(V_{i}\otimes W,W)};$$
  \item[(iii)]For any irreducible $H/\mathfrak{m}H$-module $W$, the tensor product $W\otimes W^{*}$ is a completely reducible $H/\mathfrak{m}_{\overline{\varepsilon}}H$-module.
 \end{itemize}
 \end{itemize}
\end{theorem}


We need two lemmas to prove the above theorem. The first lemma is well-known 
 (see, for example, \cite[Prop. 1.32]{MR3837537} and the discussion preceding it).
\begin{lemma}\label{lem-dimHomleq-commlcty}
Let $A$ be a $\mathbbm{k}$-algebra, $V$ a finite-dimensional irreducible $A$-module, and $W$ a finite-dimensional $A$-module.  Then 
$$[W:V]\geq \dim_{\mathbbm{k}}\text{Hom}_{A}(V,W).$$
The equality holds when $W$ is completely reducible.
\end{lemma}


\begin{lemma}
Let $H$ be a finite-dimensional Hopf algebra with the Chevalley property.\label{lem-fdchevhopf-grfus-cre} Then the Grothendieck ring $\text{Gr}(H)$ of $H$ is a fusion ring with the $\mathbb{Z}_{+}$-basis given by the set of isomorphic classes of all irreducible $H$-modules $\text{Irr}(H)=\{[V_{1}],...,[V_{m}]\}$, and the canonical regular element of $\text{Gr}(H)$ is given by $\sum_{i=1}^{m}(\dim_{\mathbbm{k}}V_{i})[V_{i}].$
\end{lemma}
\begin{proof}
By assumption, the Jacobson radical $J$ of $H$ is a Hopf ideal and $H/J$ is a semisimple Hopf algebra. 
It follows from the ring isomorphism $\text{Gr}(H)\to \text{Gr}(H/J), [V_{i}]\mapsto [V_{i}]$ that $\text{Gr}(H)$ is a fusion ring. The second conclusion follows from Lemma \ref{lem-z+ring-regel}(2).
\end{proof}
\begin{proof}
[Proof of Theorem \upshape\ref{thm-sdchevpo-geqmepi}]
(1) Let $\text{Jac}(H/\mathfrak{m}_{\overline{\varepsilon}}H)$ be the Jacobson radical of the identity fiber algebra $H/\mathfrak{m}_{\overline{\varepsilon}}H$. Denote the $H$-module $(H/\mathfrak{m}_{\overline{\varepsilon}}H)/\text{Jac}(H/\mathfrak{m}_{\overline{\varepsilon}}H)$ by $R$. Since the identity fiber algebra $H/\mathfrak{m}_{\overline{\varepsilon}}H$ has the Chevalley property, it follows from Lemma \ref{lem-fdchevhopf-grfus-cre} that
$$\mathscr{R}\coloneqq   [R]=\left[\oplus_{i=1}^{m}V_{i}^{\oplus \dim_{\mathbbm{k}}V_{i}}\right]=\sum_{i=1}^{m}(\dim_{\mathbbm{k}}V_{i})[V_{i}] \in \text{Gr}(H/\mathfrak{m}_{\overline{\varepsilon}}H)$$
is the canonical regular element of the fusion ring $\text{Gr}(H/\mathfrak{m}_{\overline{\varepsilon}}H)$. Now, Lemma \ref{lem-z+ring-regel} (1) yields
\begin{equation}\label{eq-mar2=sdmr}
\mathscr{R}^{2}=(\sum_{i=1}^{m}(\dim_{\mathbbm{k}}V_{i})^{2})\mathscr{R}=\text{Sd}(\mathfrak{m}_{\overline{\varepsilon}})\mathscr{R}.
\end{equation}
Let $\text{Irr}(H/\mathfrak{m}H)=\{[W_{1}],...,[W_{\ell}]\}$. By Corollary \ref{cor-GrirreZ+-hopf},  $\text{Gr}(H/\mathfrak{m}H)$ is an irreducible $\mathbb{Z}_{+}$-module over  $\text{Gr}(H/\mathfrak{m}_{\overline{\varepsilon}}H)$ with the fixed $\mathbb{Z}$-basis $\text{Irr}(H/\mathfrak{m}H)$. It follows that the matrix $T\in \text{M}_{\ell}(\mathbb{Z})$, defined by
$$\mathscr{R}([W_{1}],...,[W_{\ell}])=([W_{1}],...,[W_{\ell}])T,$$
has positive entries. By the Frobenius-Perron theorem, the matrix $T$ has a positive eigenvalue, and the maximal positive eigenvalue is a simple root of the characteristic polynomial of $T$. By \eqref{eq-mar2=sdmr}, $T^{2}=\text{Sd}(\mathfrak{m}_{\overline{\varepsilon}})T$. Thus, the characteristic polynomial of the matrix $T$ is given by $t^{\ell-1}(t-\text{Sd}(\mathfrak{m}_{\overline{\varepsilon}}))\in \mathbb{Z}[t]$. In particular, the trace of the matrix $T$ is $\text{Sd}(\mathfrak{m}_{\overline{\varepsilon}})$. Since the entries of the matrix $T$ are all positive integers, we see that 
$$\ell = |\text{Irr}(H/\mathfrak{m}H)|\leq \text{Sd}(\mathfrak{m}_{\overline{\varepsilon}}).$$

Note also that the trace of the matrix $T$ can be calculated by
\begin{equation}\label{eq-trT-2ndfor}
\sum_{j=1}^{\ell}[R\otimes W_{j}:W_{j}].
\end{equation}
Now, we proceed to estimate \eqref{eq-trT-2ndfor}. By Proposition \ref{prop-epiChev-mepiC}, $\text{maxSpec}C=  \text{Chev}_{\overline{\varepsilon}}(H,C)$. It follows that for each $[W_{j}]\in \text{Irr}(H/\mathfrak{m}H)$, $R\otimes W_{j}$ is a  completely reducible $H$-module. Hence, 
$$[R\otimes W_{j}:W_{j}]=\dim_{\mathbbm{k}} \text{Hom}_{H}(R\otimes W_{j},W_{j}).$$
Then
\begin{align*}
[R\otimes W_{j}:W_{j}]=&\dim_{\mathbbm{k}}\text{Hom}_{H}(R\otimes W_{j},W_{j}) \\
=&\sum_{i=1}^{m}\dim_{\mathbbm{k}}\text{Hom}_{H}(V_{i}^{\oplus \dim_{\mathbbm{k}}V_{i}}\otimes W_{j},W_{j})\\
=&\sum_{i=1}^{m}(\dim_{\mathbbm{k}}V_{i})(\dim_{\mathbbm{k}}\text{Hom}_{H}(V_{i}\otimes W_{j},W_{j})) \\
=&\sum_{i=1}^{m}(\dim_{\mathbbm{k}}V_{i})(\dim_{\mathbbm{k}}\text{Hom}_{H}(V_{i},W_{j}\otimes W_{j}^{*})) &  (\text{By } \eqref{eq-naadj-isous7})\\
\leq &\sum_{i=1}^{m}(\dim_{\mathbbm{k}}V_{i})[W_{j}\otimes W_{j}^{*}:V_{i}]  &  (\text{By Lemma }\ref{lem-dimHomleq-commlcty})\\
=&\dim_{\mathbbm{k}}(W_{j}\otimes W_{j}^{*})\\
=&(\dim_{\mathbbm{k}}W_{j})^{2}.
\end{align*}

It follows from Lemma \ref{lem-dimHomleq-commlcty} and \eqref{eq-naadj-isous7} that
$$\text{Sd}(\mathfrak{m}_{\overline{\varepsilon}})=\sum_{j=1}^{\ell}[R\otimes W_{j}:W_{j}]\leq \sum_{j=1}^{\ell}(\dim_{\mathbbm{k}}W_{j})^{2}=\text{Sd}(\mathfrak{m}).$$
\par(2) The implication (iii)$\Rightarrow$ (i) follows immediately from Lemma \ref{lem-dimHomleq-commlcty} and the proof of (1). The implication (ii) $\Rightarrow$ (iii) is obvious. It is left to prove (i) $\Rightarrow$ (ii). If $\text{Sd}(\mathfrak{m})=\text{Sd}(\mathfrak{m}_{\overline{\varepsilon}})$, then, by the proof of (1), one has
$$\dim_{\mathbbm{k}}\text{Hom}_{H}(V_{i},W_{j}\otimes W_{j}^{*})=[W_{j}\otimes W_{j}^{*}:V_{i}],\forall 1\leq i\leq m\text{ and\ } 1\leq j\leq \ell.$$
It follows that the left $H/\mathfrak{m}_{\overline{\varepsilon}}H$-module $W_{j}\otimes W_{j}^{*}$ has a submodule isomorphic to
$$\bigoplus_{i=1}^{m}V_{i}^{\oplus [W_{j}\otimes W_{j}^{*}:V_{i}]}.$$
By comparing the $\mathbbm{k}$-dimensions of these two modules immediately yields
$$W_{j}\otimes W_{j}^{*}\cong \bigoplus_{i=1}^{m}V_{i}^{\oplus [W_{j}\otimes W_{j}^{*}:V_{i}]}.$$
This shows that $W_{j}\otimes W_{j}^{*}$ is complete reducible. It follows from Lemma \ref{lem-dimHomleq-commlcty} that
$$[W_{j}\otimes W_{j}^{*}:V_{i}]=\dim_{\mathbbm{k}}\text{Hom}_{H}(V_{i},W_{j}\otimes W_{j}^{*})=\dim_{\mathbbm{k}}\text{Hom}_{H}(V_{i}\otimes W_{j},W_{j})$$
for all $1\leq j\leq \ell$.
\end{proof}
The following example shows that the two inequalities in Theorem \ref{thm-sdchevpo-geqmepi} (1) may be strict.
\begin{example}
 Let $H=H(n,t,\xi)$ be the $\mathbbm{k}$-algebra generated by $x$ and $g$ subject to the relations 
$$g^{n}=1\text{ and\ }xg=\xi gx,$$ where $n\geq 2,1\leq t\leq n$ are two coprime integers\label{eg-infiTaft}, and $\xi$ is an $n$th primitive root of unity.  
 The coalgebra structure of $H$ is defined by
$$\Delta(g)=g\otimes g,\varepsilon(g)=1\text{ and }\Delta(x)=x\otimes g^{t}+1\otimes x,\varepsilon(x)=0.$$
Then $H$ is a Hopf algebra \cite{MR2320655}, the antipode $S$ of $H$ is given by $S(g)=g^{-1}$ and $S(x)=-xg^{-t}$. Note that the coalgebra structure of $H$ here differs from that given in \cite[Ex. 2.7]{MR2320655} by a coopposite.
Sometimes $H$ is called the \textit{infinite Taft algebra}, which is  an affine prime AS-regular Hopf algebra of GK-dimension one with the center $C=\mathbbm{k}[x^{n}]$. It is direct to check that $C $ is a Hopf subalgebra of $H$ and $H$ is a free $C$-module of rank $n^{2}$. In particular, $\text{dim}_{\mathbbm{k}}H/\mathfrak{m}H=n^{2}$ for any $\mathfrak{m}\in \text{maxSpec}C$. Clearly, $H/\mathfrak{m}_{\overline{\varepsilon}}H$ is precisely the $n^{2}$-dimensional Taft algebra 
$$ H_{n}(\xi)=\mathbbm{k}\langle g,x\rangle /(g^{n}-1,xg-\xi gx, x^{n}).$$
Hence, the identity fiber algebra $H/\mathfrak{m}_{\overline{\varepsilon}}H$ is basic and $G((H/\mathfrak{m}_{\overline{\varepsilon}}H)^{\circ})$ is a cyclic group of order $n$. In  particular, $H/\mathfrak{m}_{\overline{\varepsilon}}H$ has the Chevalley property. Let $\mathfrak{m}_{\alpha}$ be the maximal ideal of $C$ generated by $(x^{n}-\alpha)$ where $\alpha\in \mathbbm{k}$. Then $\mathfrak{m}_{0}=\mathfrak{m}_{\overline{\varepsilon}}$. For $\alpha\in \mathbbm{k}^{\times}$, one has $H/\mathfrak{m}_{\alpha}H\cong \text{M}_{n}(\mathbbm{k})$. Indeed, it is direct to check that the algebra homomorphism $\Phi:H/\mathfrak{m}_{\chi_{\alpha}}H\to\text{M}_{n}(\mathbbm{k})$, given by
$$\Phi(g)=
\begin{pmatrix}
1& & & & \\
 &\xi & & &\\
& & \xi^{2} & &\\
& & &\ddots &\\
& & & &\xi^{n-1} 
\end{pmatrix},
\Phi(x)=
\begin{pmatrix}
 &\sqrt[n]{\alpha} & & & \\
 & &\sqrt[n]{\alpha} & &\\
& &  &\ddots &\\
& & & &\sqrt[n]{\alpha} \\
\sqrt[n]{\alpha}& & & &
\end{pmatrix},
$$
is bijective. It follows that for any $\alpha\in \mathbbm{k}^{\times}$, $|\text{Irr}(H/\mathfrak{m}_{\alpha}H)|=1< n$, and one has
\begin{equation}\label{Sd-intaft-max-eq}
\text{Sd}(\mathfrak{m}_{\alpha})=\begin{cases}
n, & \alpha=0,\\
n^{2}, & \alpha\neq 0.
\end{cases}
\end{equation}
Note that the identity fiber algebra $H/\mathfrak{m}_{\overline{\varepsilon}}H\cong H_{n}(\xi)$  is not semisimple, while for any $\mathfrak{m}\neq \mathfrak{m}_{\overline{\varepsilon}}\in \text{maxSpec}C$, the corresponding fiber algebra $H/\mathfrak{m}H\cong \text{M}_{n}(\mathbbm{k})$ is semisimple.
\end{example}


\subsection{Proof of Theorem \ref{intro-thmB}}\label{subsec-epChev-zerodis}
For an algebra with trace $(A,C,\text{tr})$ such that $A$ is a finitely generated free $C$-module, one defines \cite{MR1972204} the \textit{discriminant} of $A$ over $C$ with respect to the trace map $\text{tr}:A\to C$ by
\begin{equation}\label{eq-discri-def1}
\Delta(A/C;\text{tr})\coloneqq \text{det}(\text{tr}(x_{i}x_{j}))_{i,j=1}^{n}\in C,
\end{equation}
where $\{x_{1},...,x_{n}\}$ is a $C$-basis of $A$. The discriminant $\Delta(A/C;\text{tr})$ is unique up to multiplication by an invertible element in $C$. For example, when $A$ is the ring of integers $\mathcal{O}_{L}$ of an algebraic number field $L$ and $C=\mathbb{Z}$, the discriminant of $\mathcal{O}_{L}$ over $\mathbb{Z}$ with respect to the regular trace (Example \ref{eg-regulartr}) is precisely the discriminant of the algebraic number field $L$. 

In general, an algebra with trace $(A,C,\text{tr})$ does not necessarily satisfy the freeness condition. 
Instead,  one considers the \textit{discriminant ideals} \cite{MR1972204} and the \textit{modified discriminant ideals} \cite{MR3415697}.
\begin{definition}
Let $(A,C,\text{tr})$ be an algebra with trace. For any positive integer $k$,\label{def-discriid}
\begin{itemize}
\item[(1)]the \textit{$k$-discriminant ideal} $D_{k}(A/C;\text{tr})$ is the ideal of $C$ generated by the set
$$\{\text{det}(\text{tr}(a_{i}a_{j}))_{i,j=1}^{k}\mid (a_{1},...,a_{k})\in A^{k}\}.$$
\item[(2)]the \textit{modified $k$-discriminant ideal} $MD_{k}(A/C;\text{tr})$ is the ideal of $C$ generated by the set
$$\{\text{det}(\text{tr}(a_{i}b_{j}))_{i,j=1}^{k}\mid (a_{1},...,a_{k}),(b_{1},...,b_{k})\in A^{k}\}.$$
\end{itemize}
\end{definition}
Let $(A,C,\text{tr})$ be an algebra with trace. It is direct to check that
$$D_{k}(A/C;\text{tr})\subseteq MD_{k}(A/C;\text{tr})\text{ and\ }MD_{k+1}(A/C;\text{tr})\subseteq MD_{k}(A/C;\text{tr})$$
for all positive integer $k$. If $A$ can be generated by $n$ elements as a $C$-module, then 
$$MD_{k}(A/C;\text{tr})=0,$$
for all $k\geq n+1$. If, furthermore, $A$ is a free $C$-module of rank $n$, then 
\begin{equation}\label{eq-disc=ideal1}
D_{n}(A/C;\text{tr})=MD_{n}(A/C;\text{tr})=\big(\Delta(A/C;\text{tr})\big).
\end{equation}
Therefore, the discriminant can be considered as ``the top discriminant ideal'', while the remaining discriminant ideals provide information that cannot be read off from the discriminant alone. For example, see the formula \eqref{eq-CHzero-ire-sd} below.

Now, we fix an affine Cayley-Hamilton Hopf algebra $(H,C,\text{tr})$ in this subsection. We will show that if the identity fiber algebra $H/\mathfrak{m}_{\overline{\varepsilon}}H$ has the Chevalley property, then any non-empty zero locus of a discriminant ideal of $(H,C,\text{tr})$ contains the orbit of the identity element of $\text{maxSpec}C$ of the winding automorphism group action.  

In general, for an affine Cayley-Hamilton algebra $(A,C,\text{tr})$ over an algebraically closed field of characteristic zero, by \cite[Thm. 4.5 (a)]{MR1288995}, $A$ is a module-finite $C$-algebra and $C$ is an affine $\mathbbm{k}$-algebra. Recall that $\text{Sd}:\text{maxSpec}C\to \mathbb{N}$ is the square dimension function defined in \eqref{eq-Sd-function}. By \cite[Thm. 4.1 (b)]{MR3886192}, one has
\begin{equation}\label{eq-CHzero-ire-sd}
\mathcal{V}_{k}\coloneqq \mathcal{V}(D_{k}(A/C;\text{tr}))=\mathcal{V}(MD_{k}(A/C;\text{tr}))=\{\mathfrak{m}\in \text{maxSpec}C\mid \text{Sd}(\mathfrak{m})<k\},
\end{equation}
for all positive integer $k$. Therefore, it is obvious from \eqref{eq-CHzero-ire-sd} that each zero locus of a discriminant ideal of $(H,C,\text{tr})$ can be decomposed into a disjoint union of some orbits of the left winding automorphism group action. 

Now we proceed to prove the main result of this subsection.
\begin{theorem}
Let $(H,C,\text{tr})$ be an affine Cayley-Hamilton Hopf algebra such that the identity fiber algebra $H/\mathfrak{m}_{\overline{\varepsilon}}H$ has the Chevalley property. Then any non-empty zero locus of a discriminant ideal of $(H,C,\text{tr})$ contains the orbit of $\mathfrak{m}_{\overline{\varepsilon}}$ of the left winding automorphism group action.\label{thm-Chein-VDneth-cwao}
\end{theorem}

\begin{proof}
Assume that $\mathcal{V}\neq \varnothing$ is the zero locus of the $k$-discriminant ideal, that is, $\mathcal{V}=\mathcal{V}_{k}$.   Suppose $\mathfrak{m}\in \mathcal{V}_{k}$. Since $H/\mathfrak{m}_{\overline{\varepsilon}}H$ has the Chevalley property, by Theorem \ref{thm-sdchevpo-geqmepi} (1) and \eqref{eq-CHzero-ire-sd}, 
$$\text{Sd}(\mathfrak{m}_{\overline{\varepsilon}})\leq \text{Sd}(\mathfrak{m})< k.$$ 
For any point $\mathfrak{n}$ in the orbit of $\mathfrak{m}_{\overline{\varepsilon}}$ of the left winding automorphism group action, there is $\chi\in G(H^{\circ})$ such that $\mathfrak{n}=W_{l}(\chi)(\mathfrak{m}_{\overline{\varepsilon}})$. It follows that
 the winding automorphism $W_{l}(\chi):H\to H$ induces an algebra isomorphism $H/\mathfrak{n}H\cong H/\mathfrak{m}_{\overline{\varepsilon}}H$. Hence
$$
\text{Sd}(\mathfrak{n})=\text{Sd}(\mathfrak{m}_{\overline{\varepsilon}})<k.
$$
It follows from \eqref{eq-CHzero-ire-sd} that $\mathfrak{n}\in \mathcal{V}_{k}$. Hence $W_{l}(G(H^{\circ}))(\mathfrak{m}_{\overline{\varepsilon}})\subseteq \mathcal{V}_{k}$.
\end{proof}
\begin{remark}
Retain the assumptions and notation of Theorem \ref{thm-Chein-VDneth-cwao}. By Lemma \ref{lem-fdchevhopf-grfus-cre}, the Grothendieck ring $\text{Gr}(H/\mathfrak{m}_{\overline{\varepsilon}}H)$ is a fusion ring. It follows from \eqref{eq-fus-FPdim} that
$$\text{FPdim}(\text{Gr}(H/\mathfrak{m}_{\overline{\varepsilon}}H))=\text{Sd}(\mathfrak{m}_{\overline{\varepsilon}})=\sum_{[V]\in \text{Irr}(H/\mathfrak{m}_{\overline{\varepsilon}}H)}(\dim_{\mathbbm{k}}V)^{2}.$$
\end{remark}
\begin{remark}
There exist affine Cayley-Hamilton Hopf algebras $(H,C,\text{tr})$ such that the zero loci of their lowest discriminant ideals are precisely the orbits of $\mathfrak{m}_{\overline{\varepsilon}}$ under the action of the left winding automorphism group. See \cite[Thms. 5.5 (b) and 5.7 (b)]{mi2025lowest}.
\end{remark}
\begin{remark}
Some affine Cayley-Hamilton algebras are also Poisson orders in the sense of Brown-Gordon \cite{MR1989650}, including many quantum groups at roots of unity. This allows the application of Poisson geometric techniques to studying the representation theory of these algebras. A notable consequence embodying this approach is the Brown-Gordon theorem \cite[Thm. 4.2]{MR1989650}.
Furthermore, in numerous significant contexts, the affine Cayley-Hamilton algebras $(A,C,\text{tr})$ under consideration not only possess Poisson order structures, but these structures are also compatible with the trace maps $\text{tr}$. Hence, they conform to the framework of Poisson trace orders as defined by Brown-Yakimov \cite{MR4707278}.
In this setting, any zero locus of a discriminant ideal of $(A,C,\text{tr})$ is a union of symplectic cores of the Poisson variety $\text{maxSpec}C$ \cite[Cor. 2.10]{MR4707278}. For definitions and further details, see \cite{MR1989650,MR4707278}.
\end{remark}


\section{Applications}\label{sec-app}
In this section, $\mathbbm{k}$ is assumed to be an algebraically closed field of characteristic zero. We prove Theorems \ref{intro-thmC} and \ref{intro-thmD} as applications.
\subsection{Proof of Theorem \ref{intro-thmC}}\label{sec-applowest1}
We apply Theorems \ref{thm-sdchevpo-geqmepi} and \ref{thm-Chein-VDneth-cwao} to determine the level of the lowest discriminant ideal of any affine Cayley-Hamilton Hopf algebra such that its identity fiber algebra has the Chevalley property. We also give some necessary and sufficient conditions for a maximal ideal $\mathfrak{m}$ of $C$ to fall into the zero locus of the lowest discriminant ideal. 


\begin{theorem}\label{thm-lowCHHopf-minlif1}
Let $(H,C,\text{tr})$ be an affine Cayley-Hamilton Hopf algebra such that its identity fiber algebra $H/\mathfrak{m}_{\overline{\varepsilon}}H$ has the Chevalley property. Then
\begin{itemize}
\item[(1)] The lowest discriminant ideal of $(H,C,\text{tr})$ is of level 
$$\ell\coloneqq \text{FPdim}(\text{Gr}(H/\mathfrak{m}_{\overline{\varepsilon}}H))+1.$$
\item[(2)] For any $\mathfrak{m}\in \text{maxSpec}C$, the following are equivalent. 
\begin{itemize}
\item[(i)]$\mathfrak{m}\in \mathcal{V}(D_{\ell}(H/C;\text{tr}))$, where $D_{\ell}(H/C;\text{tr})$  is the lowest discriminant ideal;
\item[(ii)]$\mathfrak{m}\in \mathcal{V}(MD_{\ell}(H/C;\text{tr}))$, where $MD_{\ell}(H/C;\text{tr})$  is the lowest modified discriminant ideal;
\item[(iii)] For any irreducible $H/\mathfrak{m}H$-module $W$,  as $H$-modules
$$W\otimes W^{*}\cong \bigoplus_{i=1}^{m}V_{i}^{\oplus \dim_{\mathbbm{k}}\text{Hom}_{H}(V_{i}\otimes W,W)},$$
where $\text{Irr}(H/\mathfrak{m}_{\overline{\varepsilon}}H)=\{[V_{1}],[V_{2}],\dots,[V_{m}]\}$;
\item[(iv)]  For any irreducible $H/\mathfrak{m}H$-module $W$, $W\otimes W^{*}$ is a completely reducible $H$-module.
\end{itemize}
\end{itemize}
\end{theorem}

\begin{proof} (1) Since 
$H/\mathfrak{m}_{\overline{\varepsilon}}H$ has the Chevalley property,  it follows from Lemma \ref{lem-fdchevhopf-grfus-cre} that the Grothendieck ring $\text{Gr}(H/\mathfrak{m}_{\overline{\varepsilon}}H)$ is a fusion ring. Thus, $\text{FPdim}(\text{Gr}(H/\mathfrak{m}_{\overline{\varepsilon}}H))=\text{Sd}(\mathfrak{m}_{\overline{\varepsilon}})$. Recall that $\text{Sd}:\text{maxSpec}C\to \mathbb{N}$ is the square dimension function defined in \eqref{eq-Sd-function}.  Now, (1) follows from \eqref{eq-CHzero-ire-sd} and Theorem \ref{thm-Chein-VDneth-cwao}.

(2) This is an immediate consequence of (1), Theorem \ref{thm-sdchevpo-geqmepi} (2) and the formula \eqref{eq-CHzero-ire-sd}.
\end{proof}

\begin{remark}
Note that any basic finite-dimensional Hopf algebra has the Chevalley property. Therefore, Theorem \ref{thm-lowCHHopf-minlif1} is a generalization of \cite[Thm. B (b)(c)]{mi2025lowest}. In fact, in this case,
$$\text{FPdim}(\text{Gr}(H/\mathfrak{m}_{\overline{\varepsilon}}H))+1=|G_{0}|+1$$
where $G_{0}=G((H/\mathfrak{m}_{\overline{\varepsilon}}H)^{\circ})$ is the group of characters of the identity fiber algebra $H/\mathfrak{m}_{\overline{\varepsilon}}H$. Both the quantized coordinate rings at roots of unity and the big quantized Borel subalgebras at roots of unity fall into this situation, see \cite[Thms. 5.5 and 5.7]{mi2025lowest}.

In \cite[Def. 3.8]{mi2025lowest}, Mi-Wu-Yakimov introduced the concept of a \textit{maximally stable irreducible module}: an irreducible $H$-module $V$ is said to be maximally stable, if 
$$|\text{Stab}_{G_{0}}(V)|=(\dim_{\mathbbm{k}}V)^{2},$$
where $\text{Stab}_{G_{0}}(V)$ denotes the stabilizer of $V$ under the group action of $G_{0}$ on $\text{Irr}(H/\mathfrak{m}H)$ via tensor product. For any affine Cayley-Hamilton Hopf algebra $(H,C,\text{tr})$ whose identity fiber algebra $H/\mathfrak{m}_{\overline{\varepsilon}}H$ is basic, it was proved in \cite[Thm. 4.2 (c)]{mi2025lowest} that a maximal ideal $\mathfrak{m}$ of $C$ belongs to the zero locus of the lowest discriminant ideal of $(H,C,\text{tr})$ if and only if all irreducible $H/\mathfrak{m}H$-modules are maximally stable. This result no longer holds if the identity fiber algebra $H/\mathfrak{m}_{\overline{\varepsilon}}H$ is only assumed to have the Chevalley property. For example, consider the semisimple group algebra $H=\mathbbm{k}A_{5}$ of the alternating group $A_{5}$, and take $C=\mathbbm{k}$. 
Since $A_{5}$ is a simple group, there is only one equivalence  class of $1$-dimensional representations of $A_{5}$ over $\mathbbm{k}$. In particular, the group of characters $G_{0}$ is trivial. Now, for any irreducible $H$-module $V$ with $\dim_{\mathbbm{k}}V\geq 2$ (such an irreducible module exists as $\dim_{\mathbbm{k}}H=60$). Then $$|\text{Stab}_{G_{0}}(V)|\leq |G_{0}|=1< (\dim_{\mathbbm{k}}V)^{2}.$$
Therefore, $V$ is not maximally stable. Nevertheless, the unique maximal ideal of $C$ is of course in the zero locus of the lowest discriminant ideal of $(H,C,\text{tr}_{\text{reg}})$.
\end{remark}

Below we provide some examples. 
These example show that
the number of distinct zero loci of discriminant ideals is independent of whether the identity fiber algebra $H/\mathfrak{m}_{\overline{\varepsilon}}H$ is semisimple, or whether the global dimension of $H$ is finite. Example \ref{eg-diszero-taftinf} shows that the zero locus of a discriminant ideal of a Cayley-Hamilton Hopf algebra $(H,C,\text{tr})$ can be non-trivial.
\begin{example}\label{eg-discrzer-grpalg}
Let $G$ be a finitely generated group with a central subgroup $N$ of finite index. By Example \ref{eg-grpalg-cenex1}, for the pair $(H,C)\coloneqq (\mathbbm{k}G,\mathbbm{k}N)$, the algebra with trace $(H,C,\text{tr}_{\text{reg}})$ is a Caylay-Hamilton Hopf algebra of degree $m\coloneqq [G:N]$. Since $H/\mathfrak{m}_{\overline{\varepsilon}}H$ is semisimple, it has the Chevalley property. Then, all fiber algebras $H/\mathfrak{m}H$ are semisimple with a common $\mathbbm{k}$-dimension $m$. By \eqref{eq-CHzero-ire-sd}, 
$$\mathcal{V}(D_{k}(H/C;\text{tr}_{\text{reg}}))=\mathcal{V}(MD_{k}(H/C;\text{tr}_{\text{reg}}))=
\begin{cases} 
\varnothing, & 1\leq k\leq m, \\
\text{maxSpec}C, & k\geq m+1.
\end{cases}$$
\end{example}

\begin{example}\label{eg-16dimChe-dis-zero}
Let $(H,C)$ be the pair given in Example \ref{eg-Chevnotsembas-dim16} and retain the notation therein. Then $H$ is a $16$-dimensional Hopf algebra with the Chevalley property, and the identity fiber algebra $H/\mathfrak{m}_{\overline{\varepsilon}}H$ is a basic Hopf algebra of dimension $8$. The remaining fiber algebra $H/\mathfrak{m}H$ has a unique isomorphism class of irreducible modules with $\mathbbm{k}$-dimension $2$. 
The identity fiber algebra $H/\mathfrak{m}_{\overline{\varepsilon}}H$ has four isomorphism classes of irreducible modules. Note that $(H,C,\text{tr}_{\text{reg}})$ is a Cayley-Hamilton Hopf algebra of degree $8$. By \eqref{eq-CHzero-ire-sd},
$$\mathcal{V}(D_{k}(H/C;\text{tr}_{\text{reg}}))=\mathcal{V}(MD_{k}(H/C;\text{tr}_{\text{reg}}))=
\begin{cases} 
\varnothing, & 1\leq k\leq 4, \\
\text{maxSpec}C, & k\geq 5.
\end{cases}$$
Since $C$  is reduced (as $\mathbbm{k}$ is of characteristic zero), it follows from \eqref{eq-disc=ideal1} that $\Delta(H/C;\text{tr}_{\text{reg}})=0$.
\end{example}
\begin{example}\label{eg-diszero-taftinf}
Let $(H,C)$ be the pair given in Example \ref{eg-infiTaft} and retain the assumptions and notation therein. Then $H=H(n,t,\xi)$, where $\xi$ is a $n$th primitive  root of unity, $n\geq 2$, and the positive integer $t$ is coprime to $n$. The identity fiber algebra $H/\mathfrak{m}_{\overline{\varepsilon}}H$ is the $n^{2}$-dimensional Taft algebra $H_{n}(\xi)$ at the root of unity $\xi$. Since $H$ is a free $C$-module of rank $n^{2}$, $(H,C,\text{tr}_{\text{reg}})$ is a Cayley-Hamilton Hopf algebra of degree $n^{2}$. According to (\ref{Sd-intaft-max-eq}) and (\ref{eq-CHzero-ire-sd}), 
\begin{equation}\label{eq-zeroset-inftaft}
\mathcal{V}(D_{k}(H/C;\text{tr}_{\text{reg}}))=\mathcal{V}(MD_{k}(H/C;\text{tr}_{\text{reg}}))=
\begin{cases} 
\varnothing, & 1\leq k\leq n, \\
\{\mathfrak{m}_{\overline{\varepsilon}}\}, & n+1\leq k\leq n^{2},\\
\text{maxSpec}C, & k\geq n^{2}+1.
\end{cases}
\end{equation}
\end{example}

The zero loci of the lowest discriminant ideals of the Cayley-Hamilton Hopf algebras in Examples \ref{eg-discrzer-grpalg} and  \ref{eg-16dimChe-dis-zero} are both the entire maximal spectrum of their respective central Hopf subalgebras, whereas this is not the case for the Cayley-Hamilton Hopf algebras given by the infinite Taft algebras in Example \ref{eg-diszero-taftinf}. In the next subsection, we will see that the key reason for these two different phenomena lies in whether the Hopf algebra in question has the Chevalley property.

\subsection{Proof of Theorem \ref{intro-thmD}}\label{sec-Chev-alllrti2}
For any affine Cayley-Hamilton Hopf algebra $(H,C,\text{tr})$ with the Chevalley property, we prove that all discriminant ideals are trivial, that is, every discriminant ideal is either $0$ or $C$, which gives a necessary condition for a Cayley-Hamilton Hopf algebra to have the Chevalley property via discriminant ideals. 


\begin{theorem}\label{thm-CHHopf-Chev-dis}
Let $(H,C,\text{tr})$ be an affine Cayley-Hamilton Hopf algebra. If $H$ has the Chevalley property, then, for $\ell =\text{FPdim}(\text{Gr}(H/\mathfrak{m}_{\overline{\varepsilon}}H))+1$,
$$D_{k}(H/C;\text{tr})=MD_{k}(H/C;\text{tr})=
\begin{cases} 
C, & 1\leq k\leq \ell-1, \\
0, & k\geq \ell.
\end{cases}$$
\end{theorem}
\begin{proof}
Since $H$ has the Chevalley property, so does $H/\mathfrak{m}_{\overline{\varepsilon}}H$. Moreover, the condition (iv) in Theorem \ref{thm-lowCHHopf-minlif1} (2) is satisfied for any $\mathfrak{m}\in \text{maxSpec}C$, as the left dual of any irreducible $H$-module is still irreducible. Thus, one has 
$$\mathcal{V}(D_{k}(H/C;\text{tr}))=\mathcal{V}(MD_{k}(H/C;\text{tr}))=
\begin{cases} 
\varnothing, & 1\leq k\leq \ell-1, \\
\text{maxSpec}C, & k\geq \ell.
\end{cases}$$
Since $\text{char}\mathbbm{k}=0$,  $C$ is reduced by Cartier's Theorem. The conclusion follows.
\end{proof}
\begin{remark}
Theorem \ref{thm-CHHopf-Chev-dis}  implies that any affine Cayley-Hamilton Hopf algebra $(H,C,\text{tr})$ with the Chevalley property satisfies that
$$\varnothing =\mathcal{V}(D_{1}(H/C;\text{tr}))=\cdots=\mathcal{V}(D_{\ell-1}(H/C;\text{tr}))\subsetneq \text{maxSpec}C=\mathcal{V}(D_{\ell}(H/C;\text{tr}))=\cdots,
$$
where $\ell\coloneqq \text{FPdim}(\text{Gr}(H/\mathfrak{m}_{\overline{\varepsilon}}H))+1$. Thus, the number of distinct terms in the ascending chain of zero loci of the discriminant ideals of $(H,C,\text{tr})$ can be interpreted as a measure of how far the Hopf algebra $H$ is from having the Chevalley property.
\end{remark}

By Example \ref{eg-discrzer-grpalg}, all the zero loci of discriminant ideals of the group algebras of central extensions of finite groups by finitely generated abelian groups are trivial. This can also be regarded as a direct consequence of Theorem \ref{thm-CHHopf-Chev-dis}.

\begin{example} Let $H=H(n,t,\xi)$ be the infinite Taft algebra in Example \ref{eg-infiTaft}, where $n\geq 2, 1\leq t\leq n$ are two coprime integers, and $\xi$ is an $n$th primitive root of unity.  By (\ref{eq-zeroset-inftaft}) and Theorem \ref{thm-CHHopf-Chev-dis}, $H$ does not have the Chevalley property.
\end{example}

The following is an immediate consequence of \eqref{eq-disc=ideal1} and Theorem \ref{thm-CHHopf-Chev-dis}.
\begin{corollary}\label{cor-discri-chev-test}
Let $(H,C)$ be a pair of Hopf algebras satisfying the hypothesis (FinHopf) and $H$ is a free $C$-module.  If the discriminant $\Delta(H/C;\text{tr}_{\text{reg}})$ (with respect to the regular trace) is neither zero nor an invertible element in $C$, then $H$ does not have the Chevalley property.
\end{corollary}

In fact, for affine prime Cayley-Hamilton Hopf algebras, we can say more.
\begin{corollary}
Let $(H,C,\text{tr})$ be an affine Cayley-Hamilton Hopf algebra, and $H$ be a prime ring. Then the following conditions are equivalent. \label{cor-chev-primiff}
\begin{itemize}
\item[(1)] $H$ has the Chevalley property;
\item[(2)]All the discriminant ideals of $(H,C,\text{tr})$ are trivial;
\item[(3)]The square dimension function $\text{Sd}:\text{maxSpec}C\to\mathbb{N}$ is constant;
\item[(4)]All the fiber algebras $H/\mathfrak{m}H$ are semisimple;
\item[(5)]The identity fiber algebra $H/\mathfrak{m}_{\overline{\varepsilon}}H$ is semisimple;
\item[(6)] $H$ is commutative.
\end{itemize}
\end{corollary}
\begin{proof}
(1)$\Rightarrow$(2) This is clear by Theorem \ref{thm-CHHopf-Chev-dis}.
\par (2)$\Rightarrow$(3) It follows from (\ref{eq-CHzero-ire-sd}).
\par (3)$\Rightarrow$(4) For any $\mathfrak{m},\mathfrak{n}\in \text{maxSpec}C$, $\text{Sd}(\mathfrak{m})=\text{Sd}(\mathfrak{n})$. By \cite[Cor. III.4.10]{MR1898492}, there is a fiber algebra $H/\mathfrak{m}_{0}H$ which is semisimple, which implies that $\text{Sd}(\mathfrak{m}_{0})=\dim_{\mathbbm{k}}H/\mathfrak{m}_{0}H$. It follows from Proposition \ref{prop-equfibalg-conrp} (1) that $\text{Sd}(\mathfrak{m})=\dim_{\mathbbm{k}}H/\mathfrak{m}H$ for all $\mathfrak{m}\in \text{maxSpec}C$. Therefore $H/\mathfrak{m}H$ is semisimple for all $\mathfrak{m}\in \text{maxSpec}C$.
\par(4)$\Rightarrow$(5) Obvious.
\par(5)$\Rightarrow$(6) By Corollary \ref{cor-semisimall1}, the fiber algebra $H/\mathfrak{m}H$ is semisimple for any $\mathfrak{m}\in \text{maxSpec}C$. 
For any $\hat{\mathfrak{m}}\in \text{maxSpec}Z(H)$, $\hat{\mathfrak{m}}\cap C$ is a maximal ideal of $C$, and so $H/\hat{\mathfrak{m}}H$ is semisimple. Let $n$ be the PI-degree of $H$. By \cite[Thm. 3.1 (b)]{MR3886192}, $H/\hat{\mathfrak{m}}H\cong \text{M}_{n}(\mathbbm{k})$ as $\mathbbm{k}$-algebras. If follows that the $\mathbbm{k}$-dimension of any irreducible $H$-module is $n$. However, the trivial representation of $H$ is $1$-dimensional, which forces $n=1$. Therefore, $H$ is commutative by Posner's Theorem.
\par(6)$\Rightarrow$(1): Clear.
\end{proof}
\begin{remark}
In general, the implication (1) $\Rightarrow$ (6) in Corollary \ref{cor-chev-primiff} does not hold without the primeness condition, as Example \ref{eg-grpalg-cenex1} shows.
\end{remark}
Theorem \ref{thm-CHHopf-Chev-dis} provides a necessary condition for any affine Cayley-Hamilton Hopf algebra $(H,C,\text{tr})$ to have the Chevalley property via discriminant ideals. This is clearly not sufficient. For example, take $H$ to be the finite-dimensional Hopf algebra over $\mathbbm{k}$ in Example \ref{eg-cheloc=maxSp-notChprop} that does not have the Chevalley property, and consider $C=\mathbbm{k}$. So it is natural to ask:
\begin{que}
Is there a necessary and sufficient condition for all affine Cayley-Hamilton Hopf algebras to have the Chevalley property based on the properties of discriminant ideals?
\end{que}

\section*{Acknowledgments}
The authors thank Ruipeng Zhu for the useful discussions.
They also extend their sincere appreciation to the referees for their highly beneficial suggestions.

\end{document}